\documentclass[12pt]{article}  
\def\sq{\hbox {\rlap{$\sqcap$}$\sqcup$}}
\def\sq{\hbox {\rlap{$\sqcap$}$\sqcup$}}
\def\R{ {\rm R \kern -.31cm I \kern .15cm}}
\def\C{ {\rm C \kern -.15cm \vrule width.5pt \kern .12cm}}
\def\Z{ {\rm Z \kern -.27cm \angle \kern .02cm}}
\def\N{ {\rm N \kern -.26cm \vrule width.4pt \kern .10cm}}
\def\1{{\rm 1\mskip-4.5mu l} }
\def\lsim{\raise0.3ex\hbox{$<$\kern-0.75em\raise-1.1ex\hbox{$\sim$}}}
\def\gsim{\raise0.3ex\hbox{$>$\kern-0.75em\raise-1.1ex\hbox{$\sim$}}}
\def\noi{\noindent}

\def\beq{\begin{equation}}   \def\eeq{\end{equation}}
\def\bea{\begin{eqnarray}}  \def\eea{\end{eqnarray}}
\def\nn{\nonumber}
\def\noi{\noindent}
\def\beeq{\begin{eqnarray}} \def\eeeq{\end{eqnarray}}
\newcommand\mysection{\setcounter{equation}{0}\section}

\newcounter{hran}

\begin{document} 
\centerline{\large\bf Long range scattering for the} 
\vskip 3 truemm 
\centerline{\large\bf  Wave-Schr\"odinger system revisited} 
 \vskip 0.8 truecm

\centerline{\bf J. Ginibre}
\centerline{Laboratoire de Physique Th\'eorique\footnote{Unit\'e Mixte de
Recherche (CNRS) UMR 8627}}  \centerline{Universit\'e de Paris XI, B\^atiment
210, F-91405 ORSAY Cedex, France}
\vskip 3 truemm

\centerline{\bf G. Velo}
\centerline{Dipartimento di Fisica, Universit\`a di Bologna}  \centerline{and INFN, Sezione di
Bologna, Italy}

\vskip 1 truecm

\begin{abstract}
We reconsider the theory of scattering for the  Wave-Schr\"odinger system and more precisely the local Cauchy problem with infinite initial time, which is the main step in the construction of the wave operators. Using a method due to Nakanishi, we eliminate a loss of regularity between the Schr\"odinger asymptotic data and the Schr\"odinger solution in the treatment of that problem, in the special case of vanishing asymptotic data for the wave field.
\end{abstract}

\vskip 1 truecm
\noi 1991 MSC :  Primary 35P25. Secondary 35B40, 35Q40, 81U99.\par \vskip 2 truemm

\noi Key words : Long range scattering, Wave-Schr\"odinger system. \par 
\vskip 1 truecm

\noindent LPT Orsay 11-05\par
\noindent December 2010\par \vskip 3 truemm

\newpage
\pagestyle{plain}
\baselineskip 18pt
\mysection{Introduction}
\hspace*{\parindent} 
This paper is devoted to the theory of scattering for the Wave-Schr\"odinger (WS) system
$$\hskip 4 truecm \left\{ \begin{array}{ll}
i \partial_t u + (1/2) \Delta u = - A \ u &\hskip 4.5 truecm (1.1)\\
&\\
\sq A = |u|^2 &\hskip 4.5 truecm (1.2)
\end{array}\right .$$

\noi where $u$ and $A$ are respectively a complex valued function and a real valued function defined in space time ${I\hskip-1truemm R}^{3+1}$, 
$\Delta$ is the Laplacian in ${I\hskip-1truemm R}^3$ and $\sq = \partial_t^2 - \Delta$ is the d'Alembertian. More precisely, for arbitrarily large asymptotic data, we study the Cauchy problem at infinite initial time, which is the first step in the construction of the wave operators. This problem is complicated by the fact that the WS system is borderline long range, so that the relevant solutions of the Schr\"odinger equation contain a logarithmically diverging phase factor when $t$ tends to infinity. We previously studied that problem in \cite{3r} and we refer to the introduction of the first paper in \cite{3r} for general background. The method used in \cite{3r} is an extension of a method previously used in \cite{2r} to treat the similar case of the Hartree equation with long range potential $|x|^{-\gamma}$ with $\gamma \leq 1$. A drawback of the method used in \cite{2r}\cite{3r} is a loss of regularity of the solutions as compared with that of the asymptotic data. That defect was remedied by Nakanishi in \cite{4r} \cite{5r} for the Hartree equation in the cases $\gamma = 1$ and $1/2 < \gamma < 1$ respectively. The improvement in \cite{4r} results basically from the use of a different asymptotic parametrization of the solutions and from a clever use of the local mass conservation law for the Schr\"odinger equation. It turns out that the method used in \cite{4r}  can be extended to the WS system, unfortunately (so far) only in the special case where the wave field has zero asymptotic data. The purpose of the present paper is to present that extension, namely to solve the local Cauchy problem at infinite initial time without loss of regularity for the WS system in that special case.

In the same way as in \cite{3r}, the first step of the method consists in eliminating the wave equation by solving it for $A$ in terms of $u$. Restricting our attention to positive time and imposing the condition of vanishing asymptotic data for $A$, we obtain
$$A =  A(u,t) = - \int_t^\infty dt'\ \omega^{-1} \sin (\omega (t-t')) |u(t')|^2 \eqno(1.3)$$

\noi where $\omega = (- \Delta )^{1/2}$. We henceforth replace (1.2) by (1.3) and restrict our attention to the system (1.1) (1.3). We now introduce the relevant parametrization of $u$ needed to study the Cauchy problem at infinite time. The unitary group
$$U(t) = \exp (i(t/2)\Delta ) \eqno(1.4)$$

\noi which solves the free Schr\"odinger equation can be written as 
$$U(t) = M(t) \ D(t) \ F\ M(t) \eqno(1.5)$$

\noi where $M(t)$ is the operator of multiplication by the function
$$M(t) = \exp (ix^2/2t) \ , \eqno(1.6)$$

\noi $F$ is the Fourier transform and $D(t)$ is the dilation operator
$$D(t) = (it)^{-3/2} \ D_0(t) \eqno(1.7)$$

\noi where
$$\left ( D_0 (t) f \right ) (x) = f(x/t)\ . \eqno(1.8)$$

\noi For any function $g$ of space time, we define
$$\widetilde{g}(t) = U(t)^* \ g(t) \ . \eqno(1.9)$$

\noi We first perform a pseudoconformal inversion on $u$, namely
$$u(t) = M(t) \ D(t) \ \overline{u_c} (1/t) \eqno(1.10)$$

\noi or equivalently
$$\widetilde{u}(t) = \overline{F \widetilde{u}_c (1/t)}\ , \eqno(1.11)$$

\noi thereby replacing the Cauchy problem at infinite initial time for $u$ by the Cauchy problem at initial time zero for $u_c$. Correspondingly we replace $A(t)$ by $B(t)$ defined by
$$A(t) = t^{-1} \ D_0(t)\ B(1/t) \ . \eqno(1.12)$$

\noi The system (1.1) (1.3) is then replaced by
$$i \partial_t \ u_c = - (1/2) \Delta u_c - t^{-1} \ B(u_c) u_c \ , \eqno(1.13)$$

$$B(u_c,t) = \int_1^\infty d\nu \ \nu^{-3} \omega^{-1} \sin (\omega (\nu - 1)) D_0(\nu ) |u_c(t/\nu)|^2 \ . \eqno(1.14)$$

\noi The variables $u_c$ and/or $t$ in $B$ will be partly omitted when no confusion can arise, as for instance in (1.13). We now parametrize $u_c$ in terms of an amplitude $v$ and a phase $\varphi$ according to
$$\widetilde{u}_c(t) = \exp (- i \varphi (t)) v(t) \eqno(1.15)$$

\noi or equivalently
$$u_c(t) = U_\varphi (t)\ v(t) \eqno(1.16)$$

\noi where
$$U_\varphi (t) = U(t) \exp (- i \varphi (t))\ . \eqno(1.17)$$

\noi That parametrization is the same as that used in \cite{4r} and differs from that used in \cite{3r} where the phase factor was introduced in $u_c$ instead of $\widetilde{u}_c$. The equation (1.13) then becomes the following equation for $v$~: 
$$i \partial_t v = - \left ( t^{-1} U_\varphi^* \ B(u_c) U_\varphi + \partial_t \varphi \right ) v \ . \eqno(1.18)$$

The role of the phase $\varphi$ is to cancel the singularity at $t=0$ of the last term in (1.13), so that (1.18) can be solved with $v$ continuous at $t =0$ and with initial condition $v(0) = v_0$. It will then turn out that $B(u_c)$ tends to $B (v_0)$ when $t$ tends to zero. The cancellation will be ensured by imposing
$$\partial_t \varphi = - t^{-1} B(v_0) \eqno(1.19)$$

\noi which together with the (arbitrary) initial condition $\varphi (1) = 0$, yields 
$$\varphi (t) = - (\ell n \ t) B(v_0) \ . \eqno(1.20)$$

\noi The equation for $v$ then becomes
$$i \partial_t v = - t^{-1} \left ( U_\varphi^* \ B(u_c) U_\varphi - B(v_0)\right ) v \equiv t^{-1} \ L(v) v \eqno(1.21)$$

\noi with
$$L(v) = -  \left ( U_\varphi^* \ B(u_c) U_\varphi - B(v_0)\right )  \ . \eqno(1.22)$$

\noi We shall also need the partially linearized equation for $v'$
$$i \partial_t v' =  t^{-1} \ L(v)v'\ .  \eqno(1.23)$$

The method consists in first solving the Cauchy problem with initial time zero for the linearized equation (1.23). One then shows that the map $v \to v'$ thereby defined is a contraction in a suitable space in a sufficiently small time interval. This solves the Cauchy problem with initial time zero for the nonlinear equation (1.21). One then translates the results through the change of variables (1.15) to solve the Cauchy problem with initial time zero for the system (1.13) (1.14) or equivalently with infinite initial time for the system (1.1) (1.3). The final result can be stated as the following proposition, which is a slightly shortened rewriting of Propositions 6.1-6.3. We need the notation
$$FH^\rho = \left \{ u \in {\cal S}' : F^{-1} u \in H^\rho \right \}\ .$$

\vskip 3 truemm
\noi {\bf Proposition 1.1.} {\it Let $1 < \rho < 3/2$. \par

(1) Let $u_0\in FH^\rho$ and define
$$\varphi (t) = - \ell n \ t\ B (\overline{Fu_0})\ . \eqno(1.24)$$

\noi Then there exists $T_\infty > 0$ and there exists a unique solution $u$ of the system (1.1) (1.3) such that $\widetilde{u} \in {\cal C}([T_\infty , \infty ), FH^\rho)$ and such that $w$ defined by
$$w(t) = F^{-1} \exp (-i \varphi (1/t)) F \widetilde{u}(t) \eqno(1.25)$$

\noi satisfy
$$w(t) \to u_0 \quad {\rm in} \quad FH^\rho \quad {\rm when }\quad t \to \infty\  . \eqno(1.26)$$

\noi Furthermore $w \in {\cal C} ([T_\infty , \infty ), FH^\rho )$, the map $u_0 \to w$ is continuous from $FH^\rho$ to\break\noindent $L^\infty ([T_\infty , \infty ), FH^\rho)$ and $u$ satisfies the estimate 
$$\parallel \widetilde{u} (t) ; FH^\rho \parallel \ \leq a_1 (1 + |\ell n \ t|)^2$$

\noi for some $a_1 \geq 0$ and for all $t \geq T_\infty$.\\

(2) Let $T_\infty > 0$. Let $u$ be a solution of the system (1.1) (1.3) such that $\widetilde{u} \in {\cal C} ([T_\infty , \infty ), FH^\rho)$ and that $u$ satisfy the estimate 
$$\parallel \widetilde{u} (t) ; FH^\rho \parallel \ \leq a_1 (1 + |\ell n \ t|)^\alpha$$

\noi for some $a_1, \alpha \geq 0$ and for all $t \geq T_\infty$. Then there exists $u_0 \in FH^\rho$ such that $w$ defined by (1.25) (1.24) satisfies (1.26). Furthermore $w  \in {\cal C} ([T_\infty , \infty ), FH^\rho)$.\\

(3) Let $T_\infty > 0$. Let $u_i$, $i = 1,2$ satisfy the assumptions of Part (2) and assume that $u_1(t) -u_2(t) \to 0$ in $L^2$ when $t \to \infty$. Then $u_1 = u_2$.}\\

Part (2) of Proposition 1.1 is the converse of Part (1) in the sense that all sufficiently regular solutions of the system (1.1) (1.3) can be recovered by the construction of Part (1). As an application of that result, one obtains in Part (3) a uniqueness result for that system not making any reference to the parametrization (1.15). \par

The lower bound $\rho > 1$ in Proposition 1.1 is essential for the success of the method. Actually the value $\rho =1$ is critical in a natural sense. On the other hand the upper bound $\rho < 3/2$ is imposed for convenience only and could be dispensed with at the expense of complicating the estimates. \par

We then briefly comment on two questions which we do not consider in this paper. Firstly we do not extend the solutions of the local Cauchy problem to global ones, since the norms in Proposition 1.1 are ill adapted to the wave equation and therefore do not readily allow for globalisation. Secondly we do not consider the converse question of solving the Cauchy problem for (1.1) (1.3) up to infinity in time, starting from a (sufficiently large) finite initial time. The reason is that this problem is no longer a Cauchy problem for the original system (1.1) (1.2) after the latter has been reduced to (1.1) (1.3) by imposing that the wave field vanishes at infinity in time. \par

This paper is organized as follows. In Section 2, we introduce some notation and we collect a number of estimates which are used throughout this paper. In Section 3, we study the Cauchy problem for the linearized equation (1.23) with initial time $t_0 \geq 0$. In Section 4, we solve the Cauchy problem with initial time zero for the nonlinear equation (1.21). In Section 5, we prove the continuity of the solutions of (1.21) with respect to the initial data. In Section 6 we reformulate the previous results as results on the Cauchy problem for the system (1.13) (1.14) and we derive an additional uniqueness result for that system not making any reference to $v$. 

\mysection{Notation and preliminary estimates} 
\hspace*{\parindent}
In this section we introduce some notation and we collect a number of
estimates which will be used throughout this paper. We denote by
$\parallel \cdot \parallel_r$ the norm in $L^r \equiv
L^r({I\hskip-1truemm R}^{3})$.  For any interval $I$ and any Banach space $X$ we denote by ${\cal C}(I,
X)$ (resp. ${\cal C}_w(I,X)$, ${\cal C}^L(I,X))$ the space of strongly (resp. weakly, Lipschitz) continuous functions from $I$ to $X$ and by
$L^{\infty} (I, X)$ the space of measurable essentially bounded
functions from $I$ to $X$. For real numbers $a$ and $b$ we use the
notation $a \vee b = {\rm Max}(a,b)$ and $a\wedge b = {\rm Min} (a,b)$. \par

We shall use the Sobolev spaces $\dot{H}_r^\sigma$ and $H_r^\sigma$ defined for
$- \infty < \sigma < + \infty$, $1 \leq r \leq \infty$ by
$$\dot{H}_r^\sigma = \left \{ u:\parallel u;\dot{H}_r^\sigma\parallel \ \equiv \ \parallel \omega^\sigma u\parallel_r \ <
\infty \right \}$$

\noi and

$$H_r^\sigma = \left \{ u:\parallel u;H_r^\sigma\parallel \ \equiv \ \parallel <\omega>^\sigma u\parallel_r \ <
\infty \right \}$$

\noi where $\omega = (- \Delta)^{1/2}$ and $< \cdot > = (1 + |\cdot |^2)^{1/2}$. The subscript $r$ will
be omitted if $r = 2$ and we shall use
the notation 
$$\parallel \omega^{\sigma \pm 0} u \parallel_2\ = \left ( \parallel \omega^{\sigma + \varepsilon} u \parallel_2\ \parallel \omega^{\sigma - \varepsilon} u \parallel_2\right )^{1/2} \quad \hbox{for some $\varepsilon > 0$\ .}$$  

We shall use extensively the following Sobolev inequalities, stated
here in ${I\hskip-1truemm R}^n$, but to be used only for $n = 3$. \\

\noi {\bf Lemma 2.1.} {\it Let $1 < q, r < \infty$, $1 < p \leq \infty$ and $0 \leq \sigma < \rho$. If $p = \infty$, assume that $\rho - \sigma > n/r$. Let $\theta$ satisfy $\sigma /\rho \leq \theta \leq 1$ and
$$n/p - \sigma = (1 - \theta )n/q + \theta (n/r - \rho ) \ .$$

\noi Then the following inequality holds}
\beq  
\label{2.1e}
\parallel \omega^\sigma u\parallel_p \ \leq C  \parallel u \parallel_q^{1 - \theta} 
\ \parallel \omega^\rho u \parallel_r^\theta \ .
\eeq

We shall also use extensively the following Leibnitz and commutator estimates. \\

\noi {\bf Lemma 2.2.} {\it Let $1 < r,r_1,r_3 < \infty$ and
$$1/r = 1/r_1 + 1/r_2 = 1/r_3 + 1/r_4 \ .$$

\noi Then the following estimates hold for $\sigma \geq 0$~:}
\beq  
\label{2.2e}
\parallel \omega^\sigma (uv) \parallel_r \ \leq C \left ( \parallel  \omega^\sigma u \parallel_{r_1} 
\ \parallel v \parallel_{r_2} + \parallel  \omega^\sigma v \parallel_{r_3} \  
\parallel u \parallel_{r_4} \right )
\ .  \eeq 

An easy consequence of Lemmas 2.1 and 2.2 is the following estimate.\\

\noi {\bf Lemma 2.3.} {\it Let $1 < r < \infty$, $\sigma_1 \geq 0$,  $\sigma_2 \geq 0$, $\sigma_1 + \sigma_2 = \sigma + n - n/r$, such that $-n + n/r < \sigma \leq \sigma_1 \wedge \sigma_2 \leq \sigma_1 \vee \sigma_2 < n/2$. Then} 
\beq  
\label{2.3e}
\parallel \omega^\sigma (uv) \parallel_r \ \leq C  \parallel  \omega^{\sigma_1} u \parallel_2
\ \parallel  \omega^{\sigma_2} v \parallel_2  \  .
\eeq
\vskip 5 truemm

\noi For $0 \leq \sigma < n/2$, we define the space $M^\sigma = L^\infty \cap \dot{H}_{n/\sigma}^\sigma$ with norm
\beq  
\label{2.4e}
\parallel f;M^\sigma \parallel  \ = \  \parallel f \parallel_\infty \   + \ \parallel \omega^\sigma f  \parallel_{n/\sigma} \  .
 \eeq

\noi By Lemmas 2.1 and 2.2

$$\parallel f;M^\sigma \parallel  \ \leq \  \parallel f \parallel_\infty \   + \ C \parallel \omega^{n/2} f  \parallel_2$$
\beq
\label{2.5e}
\leq C \left ( \parallel  \omega^{n/2 + \varepsilon} f  \parallel_2 \ \parallel \omega^{n/2 - \varepsilon} f \parallel_2 \right )^{1/2} \equiv \ C\parallel \omega^{n/2\pm 0} f \parallel_2
\eeq

\noi and for $0 \leq \sigma ' \leq \sigma < n/2$
\beq
\label{2.6e}
\parallel f u ;\dot{H}^{\sigma '} \parallel  \ \leq C  \parallel f; M^\sigma \parallel \    \ \parallel u; \dot{H}^{\sigma '}   \parallel \  .
\eeq

We shall also need some commutator estimates.\\

\noi {\bf Lemma 2.4.} {\it Let $\lambda \geq 0$, let $\sigma_1 , \sigma_2$ satisfy $0 \leq \sigma_1 , \sigma_2 < n/2$ and\break\noindent  $( \lambda - 1) \vee 0 < \sigma_1 + \sigma_2 \leq \lambda + n/2$. 

Then the following estimate holds
\beq  
\label{2.7e}
|< u, [\omega^\lambda , f] v>| \leq \ C \parallel  \omega^\beta f \parallel_2\   \parallel  \omega^{\sigma_1} u \parallel_2
\ \parallel  \omega^{\sigma_2} v \parallel_2  \  ,
\eeq
\noi where $\beta = \lambda + n/2 - \sigma_1 - \sigma_2$.}\\

\noi Lemma 2.4 is a variant of Lemma 3.6 in \cite{4r} and is proved by a minor variation of the proof of the latter. \par

We now give some estimates of $B$ defined by (1.14). From now on we take $n=3$. \\

\noi {\bf Lemma 2.5.} {\it Let $0 < \sigma < 3/2$ and $u \in {\cal C}((0, 1] , \dot{H}^\sigma )$. Then 
\beq  
\label{2.8e}
\parallel \omega^{2\sigma - 1/2} B(u,t) \parallel_2 \ \leq C  \int_1^\infty d\nu \ \nu^{-2\sigma}  \ \parallel  \omega^{\sigma} u (t/ \nu )\parallel_2^2  \  .
\eeq

\noi Let $2 \leq r \leq 7$. Then 
\beq
\label{2.9e}
\parallel B(u,t) \parallel_r  \ \leq C \int_1^\infty d\nu \ \nu^{-1 + 1/r} (\nu - 1)^{-1+2/r} \ \parallel  \omega^{\sigma} u (t/\nu )\parallel_2^2  
\eeq

\noi with $\sigma = 1/2 - 1/2r$.} \\

\noi {\bf Proof.} From (1.14) and from the identity
\beq
\label{2.10e}
\parallel \omega^\alpha D_0 (\nu ) f \parallel_r \  = \nu^{-\alpha + n/r} \ \parallel \omega^\alpha f \parallel_r
\eeq

\noi we obtain
\beq  
\label{2.11e}
\parallel \omega^{2\sigma - 1/2} B(u,t) \parallel_2 \ \leq C  \ \int_1^\infty d\nu \ \nu^{-2\sigma}  \ \parallel  \omega^{2\sigma - 3/2} |u (t/ \nu )|^2\parallel_2  
\eeq

\noi from which (\ref{2.8e}) follows by Lemma 2.3.\par

From the well known dispersive estimates for the wave equation \cite{6r}, we obtain
\beq
\label{2.12e}
\parallel B(u, t) \parallel_r  \ \leq C \int_1^\infty d\nu \ \nu^{-3} (\nu - 1)^{-1+2/r} \ \parallel  \omega^{1-4/r} D_0(\nu ) |u(t/\nu )|^2\parallel_{\bar{r}} \  .
\eeq

\noi for $2 \leq r \leq \infty$, where $1/r + 1/\overline{r} = 1$, from which (\ref{2.9e}) follows by (\ref{2.10e}) and Lemma 2.3. \par 
\hfill $\sq$

\noi {\bf Remark 2.1.} Lemma 2.5 says nothing on the convergence of the integrals over $\nu$, which may well be infinite. In the applications, $u(t)$ will be bounded in time up to logarithmic factors near $t=0$, so that convergence will be ensured for $\sigma > 1/2$ in (\ref{2.8e}) and for $r > 3$ in (\ref{2.9e}).\par

In order to obtain further estimates on $B(u)$, we need additional assumptions on $u$, namely the fact that $u$ satisfies some linear Schr\"odinger equation. The following estimate, which holds for any space dimension $n \geq 2$, plays an essential role in \cite{4r} and in this paper.\\

\noi {\bf Lemma 2.6.} {\it Let $n \geq 2$, let $1/2 < \sigma < n/2$, let $I$ be an interval and let $u \in {\cal C}(I, H^\sigma )$ be a solution of the equation
\beq
\label{2.13e}
i \partial_t u + (1/2) \Delta u = Vu
\eeq

\noi in I for some real $V \in L_{loc}^\infty (I, L^\infty )$. Then for any $t_1, t \in I$, $t_1 \leq t$, the following estimate holds~:}
\beq
\label{2.14e}
\parallel \omega^{2\sigma - 2- n/2} \left ( |u(t)|^2 - |u(t_1)|^2\right )  \parallel_2 \ \leq C    \int_{t_1}^t dt' \parallel \omega^\sigma u(t') \parallel_2^2  \  .
\eeq

\noi {\bf Sketch of proof.} The formal local conservation law
\beq
\label{2.15e}
\partial_t  | u|^2 = - \ {\rm Im} \ \overline{u} \Delta u 
\eeq

\noi implies
$$ \partial_t \ < | u|^2, \psi >\  = - i/2 \ <u, [ \Delta , \psi ] u>$$

\noi for any test function $\psi$ of the space variable. Integrating over time and estimating the right hand side by Lemma 2.4 with $\lambda = 2$, $\sigma_1 = \sigma_2 = \sigma$ yields
$$|<|u(t)|^2 - |u(t_1)|^2, \psi >|  \leq C \int_{t_1}^t dt' \parallel \omega^\sigma u(t') \parallel_2^2 \ \parallel \omega^{n/2 + 2 - 2\sigma} \psi \parallel_2$$

\noi from which (\ref{2.14e}) follows by duality. The formal proof can be made rigorous under the regularity assumptions made on $u$ and $V$.\par
\hfill $\sq$

We now exploit Lemma 2.6 to derive another estimate on $B$. \\

\noi {\bf Lemma 2.7.} {\it Let $1/2< \sigma < 3/2$, let $I = (0,T]$ and let $u \in {\cal C}(I, H^\sigma )$ be a solution of the equation (\ref{2.13e}) in $I$ for some real $V \in L_{loc}^\infty (I, L^\infty )$. Then for $0 < t_1 \leq t \leq T$, the following estimate holds~:} 
\beq  
\label{2.16e}
\parallel \omega^{2\sigma - 5/2} \left ( B(u,t) - B(u, t_1)\right ) \parallel_2 \ \leq C  (3 - 2\sigma)^{-1} t \  \int_1^\infty d\nu \ \nu^{1-2\sigma}  \ \parallel  \omega^{\sigma} u (t/\nu )\parallel_2^2  \  .
\eeq

\noi {\bf Proof.} From (1.14) and (\ref{2.10e}) we estimate
$$\parallel \omega^\beta  \left ( B(u,t) - B(u, t_1)\right ) \parallel_2 \ \leq  \int_1^\infty d\nu \ \nu^{-1/2-\beta}  \ \parallel  \omega^{\beta- 1} \left ( |u (t/\nu )|^2 - |u(t_1/\nu )|^2\right ) \parallel_2$$

\noi which by (\ref{2.14e}) is continued as
\beq
\label{2.17e}
\cdots \leq C \int_1^\infty d\nu \ \nu^{2-2\sigma} \int_{t_1/\nu}^{t/\nu} dt'  \parallel \omega^\sigma u(t') \parallel_2^2
\eeq

\noi with $\beta = 2\sigma - 5/2$. Changing variables from $(t', \nu )$ to $(s, \nu ')$ defined by $t' = t/\nu '$, $\nu = \nu 's$ yields
$$\cdots = C \ t \int_1^\infty d\nu ' \ \nu{'}^{1 - 2\sigma} \parallel \omega^\sigma u(t/\nu ') \parallel_2^2 \int_{(1/\nu ' ) \wedge (t_1/t)}^1 ds\ s^{2-2\sigma}$$
$$= C\ t  \int_1^\infty d\nu  \ \nu^{1 - 2\sigma} \parallel \omega^\sigma u(t/\nu ) \parallel_2^2 (3 - 2 \sigma )^{-1} \left ( 1 - \left ( (1/\nu ) \wedge (t_1/t) \right )^{3 - 2\sigma}\right )$$

\noi which implies (\ref{2.16e}).\par\nobreak
\hfill $\sq$

The great advantage of Lemma 2.7 over Lemma 2.5 is the fact that the RHS of (\ref{2.16e}) tends to zero when $t_1, t \to 0$ under suitable assumptions on $u$. In order to allow for more flexibility, we interpolate between (\ref{2.8e}) and (\ref{2.16e}), as stated in the following lemma.\\

\noi {\bf Lemma 2.8.} {\it Under the assumptions of Lemma 2.7, the following estimate holds for $0 < t_1 \leq t \leq T$ and $0 \leq \theta \leq 1$~:}
$$\parallel \omega^{2\sigma - 1/2-2\theta} \left ( B(u,t) - B(u, t_1)\right ) \parallel_2 \ \leq C  (3 - 2\sigma)^{-\theta} t^{\theta}$$
\beq
\label{2.18e}
\int_1^\infty d\nu \ \nu^{1-2\sigma} \left (  \parallel  \omega^{\sigma} u (t/\nu )\parallel_2^2\ +  \parallel  \omega^{\sigma} u (t_1/\nu )\parallel_2^2\right )  \  .
\eeq

\noi {\bf Proof.} Interpolating between (\ref{2.8e}) and (\ref{2.16e}) yields 
$$\parallel \omega^{2\sigma - 1/2 - 2\theta} \left ( B(u,t) - B(u, t_1)\right ) \parallel_2 \ \leq C  (3 - 2\sigma)^{-\theta} t ^{\theta}$$
$$\left \{  \int_1^\infty d\nu \ \nu^{1-2\sigma}  \ f(t/\nu ) \right \}^\theta \left \{ \int_1^\infty  d\nu \ \nu^{-2\sigma} \left ( f(t/\nu ) + f(t_1 /\nu )\right ) \right \}^{1 - \theta}$$ 

\noi with $f(t) = \ \parallel \omega^{\sigma} u (t)\parallel_2^2$, from which (\ref{2.18e}) follows.\par\nobreak 
\hfill $\sq$

In order to handle the phase factor occurring in (1.15), we shall need some phase estimates. The basic estimate is best expressed in terms of homogeneous Besov spaces \cite{1r}. The following lemma is a variant of Lemma 3.3 in \cite{4r}.\\

\noi {\bf Lemma 2.9.} {\it Let $\varphi$ be a real function. Let $\sigma > 0$ and $1 \leq q, r < \infty$. Then the following estimate holds~:
\beq
\label{2.19e}
\parallel \left ( \exp (i \varphi ) - 1 \right ) ; \dot{B}_{r, q}^\sigma \parallel \ \leq \ C \parallel \varphi ; \dot{B}_{r,q}^\sigma \parallel \left ( 1 \ + \ \parallel \varphi ; \dot{B}_{\infty , \infty}^0 \parallel \right )^{[\sigma ]}
\eeq

\noi where $[\sigma ]$ is the integral part of $\sigma$.}\\

We shall use extensively the following special case of (\ref{2.19e})
\beq  
\label{2.20e}
\parallel \omega^{\sigma} (\exp (i \varphi ) - 1)  \parallel_2 \ \leq C  \parallel  \omega^{\sigma} \varphi \parallel_2 \left ( 1 \ + \  \parallel  \omega^{n/2} \varphi \parallel_2  \right )^{[\sigma ]}\   .
\eeq

\noi Lemma 2.9 implies the following estimate of the $M^\sigma$ norm of $\exp (i \varphi ) - 1$ defined by (\ref{2.4e}) for $0 < \sigma < n/2$~: 
\bea  
\label{2.21e}
&&\parallel  (\exp (i \varphi ) - 1 ); M^\sigma  \parallel \ \leq C  \left ( 1 \ \wedge \parallel  \varphi \parallel_\infty\ + \ \parallel  \varphi ; \dot{B}_{n/\sigma , 2}^{\sigma} \parallel \left ( 1 \ + \  \parallel \varphi ; \dot{B}_{\infty , \infty}^0 \parallel  \right )^{[\sigma ]}\right )   \nn \\
&&\leq \ C \left ( 1 \ \wedge \ \parallel \varphi \parallel_{\infty} \ + \ \parallel \omega^{n/2} \ \varphi \parallel_2 \left ( 1 \ + \  \parallel \omega^{n/2} \ \varphi \parallel_2 \right )^{[\sigma ]}\right ) \ .
\eea													

\noi We now exploit the previous abstract lemmas to derive some estimates of $B(u_c, t)$ defined by (1.14) with $u_c$ expressed by (1.16) and (1.20) with $v_0 = v(0)$.\\

\noi {\bf Lemma 2.10.} {\it Let $1 < \sigma < 3/2$, let $I = (0,T]$ and let $v \in L^{\infty} (I, H^\sigma ) \cap {\cal C} ([0, T], L^2)$. Let $u_c$ be defined by (1.16) (1.20) with $v_0 = v(0)$ and satisfy the equation
\beq
\label{2.22e}
i \partial_t u_c + (1/2) \Delta u_c = Vu_c
\eeq

\noi for some real $V \in L_{loc}^\infty (I, L^\infty )$. Let $0 < \theta \leq 1$. Then $B(u_c, t)$ defined by (1.14) satisfies the estimate
\bea
\label{2.23e}
&&\parallel \omega^{2\sigma - 1/2 - 2 \theta} \left ( B(u_c,t) - B(u_c, t_1) \right )  \parallel_2 \ \leq C  (\sigma , \theta )t^{\theta}\nn \\
&&\times \left ( 1 \ + \  \parallel \nabla v_0 \parallel_2^2 \left ( 1 + |\ell n \ t| \right ) \right )^4 \ \parallel v ; L^\infty \left ( (0, t], \dot{H}^\sigma \right ) \parallel^2
\eea

\noi for $0 \leq t_1 \leq t \leq T$ and
\beq
\label{2.24e}
B(u_c , 0) \ \equiv \lim_{t\to 0} B (u_c, t) = B(v_0)\  .  
\eeq

The limit in (\ref{2.24e}) holds in all the norms appearing in (\ref{2.23e}), in the sense that $B(u_c, t) - B(v_0)$ tends to zero when $t \to 0$ in $\dot{H}^\beta$ for $-1/2 <  \beta < 2\sigma - 1/2$.}\\

\noi {\bf Remark 2.2.} By Lemma 2.5, $B(v_0) \in \dot{H}^{1/2+ 0} \cap \dot{H}^{2\sigma - 1/2}$. From that fact and from the convergence of $B(u_c)$ to $B(v_0)$, it follows that $B(u_c) \in L^\infty (I, \dot{H}^{1/2 + 0} \cap \dot{H}^{2\sigma - 1/2 - 0})$. Note however that the limit in (\ref{2.24e}) holds in some norms which are not expected to be finite for $B(v_0)$, typically the $\dot{H}^\beta$ norm for $\beta \leq 1/2$.\\

\noi {\bf Proof.} We first estimate
\bea
\label{2.25e}
\parallel \omega^{\sigma} u_c(t)\parallel_2 &\leq& \ \parallel e^{i\varphi}; M^\sigma \parallel \ \parallel \omega^\sigma v(t) \parallel_2 \nn \\
&\leq& \left ( 1 + a_0^2 |\ell n\ t|\right )^2 \ \parallel \omega^\sigma v(t)\parallel_2 
\eea

\noi with $a_0 = \parallel \nabla v_0 \parallel_2$, by (\ref{2.21e}) and Lemma 2.5. Substituting (\ref{2.25e}) into (\ref{2.18e}) yields
$$\parallel \omega^\beta  \left ( B(t) - B(t_1)\right ) \parallel_2 \ \leq C  (3 - 2\sigma)^{-\theta} t^{\theta}\parallel v;L^\infty ((0, t], \dot{H}^\sigma )\parallel N$$

\noi where $\beta = 2\sigma - 1/2 - 2\theta$ and			
\bea
\label{2.26e}
N &=& \int_1^\infty d\nu \ \nu^{1-2\sigma} \left ( 1 + a_0^2 \left ( |\ell n \ t| + |\ell n \ t_1| + \ell n \ \nu \ \right ) \right )^4 \nn\\
&\leq&C(\sigma ) \left ( 1 + a_0^2 \left ( 1 + |\ell n\ t|  + |\ell n \ t_1|\right ) \right )^4 
\eea

\noi where 
$$C (\sigma ) = C \int_1^\infty d\nu \ \nu^{1-2\sigma}(1 + \ell n \ \nu)^4\ .$$

We next estimate
$$\parallel \omega^\beta  \left ( B(t) - B(t_1)\right ) \parallel_2 \ \leq \sum_{0 \leq j \leq \ell - 1} \parallel \omega^\beta  \left ( B\left ( t2^{-j}\right ) - B\left ( t 2^{-(j+1)} \right  )\right ) \parallel_2$$
$$+ \  \parallel \omega^\beta  \left ( B\left ( t2^{-\ell}\right ) - B\left ( t_1 \right ) \right ) \parallel_2$$

\noi for $2^{- (\ell + 1)} \leq t_1/t \leq 2^{-\ell}$. We estimate each term in the sum by (\ref{2.26e}) with $(t, t_1)$ replaced by $(t2^{-j}, t2^{-(j+1)})$ or by $(t2^{-\ell}, t_1)$ thereby obtaining
$$\parallel \omega^\beta  \left ( B(t) - B(t_1)\right ) \parallel_2 \ \leq C( \sigma ) (3 - 2\sigma )^{-\theta} \ t^\theta \parallel v; L^\infty ((0, t], \dot{H}^\sigma ) \parallel$$
$$\times \sum_{j \geq 0} 2^{-j\theta } \left ( 1 + a_0^2 (1 + |\ell n\ t| + j \ell n\ 2)\right )^4$$

\noi from which (\ref{2.23e}) follows for $0 < t_1 \leq t$ with
$$C(\sigma , \theta ) = C(\sigma ) (3 - 2\sigma )^{-\theta} \sum_{j \geq 0} 2^{-j \theta}\ j^4 \ .$$

\noi Note that $C (\sigma , \theta )$ blows up when $\sigma \to 1$ or $\sigma \to 3/2$ or $\theta \to 0$.

We next prove that $B (u_c, t)$ tends to $B(v_0)$ in $\dot{H}^1$. In fact 
\beq
\label{2.27e}
\parallel \omega \left ( B(t) - B(v_0)\right ) \parallel_2 \ \leq \int d\nu \ \nu^{-3/2} \parallel |u_c(t/ \nu )|^2 - |v_0|^2 \parallel_2
\eeq

\noi by (\ref{2.11e}) with $\sigma = 3/4$. Now
$$|u_c(t/\nu )|^2 - |v_0|^2 = {\rm Re} \left ( (U^* - 1) e^{i \varphi} \ \overline{v}\right ) (U + 1) e^{-i \varphi} \ v + \ {\rm Re} ( \overline{v} - \overline{v}_0) (v + v_0)$$

\noi so that
$$\parallel |u_c (t/\nu )|^2 - |v_0|^2 \parallel_2 \ \leq \ \parallel (U^* - 1) e^{i\varphi} \  \overline{v} \parallel_3 \ \parallel (U+1) e^{-i \varphi} \ v \parallel_6$$
$$+ \ \parallel v - v_0 \parallel_2^{1/2}\ \parallel v - v_0 \parallel_6^{1/2} \ \parallel v + v_0 \parallel_6$$
$$\leq C \left \{ (t/\nu)^{1/4} \parallel \nabla e^{-i\varphi} \ v\parallel_2^2 \ + \ \parallel v - v_0 \parallel_2^{1/2} \left ( \parallel \nabla v \parallel_2\ + \ \parallel \nabla v_0 \parallel_2 \right )^{3/2} \right \}$$
$$\leq C \left \{ (t/\nu )^{1/4} \left ( 1 + |\ell n\ t) \parallel \nabla B_0 \parallel_3 \right )^2 \ \parallel \nabla v \parallel_2^2 + \cdots \right \}$$

\noi which tends to zero when $t \to 0$ for $v\in L^\infty (I, H^1)$ and $v(t)$ tending to $v_0$ in $L^2$. Together with the estimate (\ref{2.23e}) for $0 < t_1 \leq t \leq T$, this proves that the same estimate also holds for $t_1 = 0$ with $B(u_c, 0) = B(v_0)$ by an appropriate abstract argument. \par
\hfill $\sq$

\mysection{The linearized Cauchy problem for v}
\hspace*{\parindent} 
In this section we study the Cauchy problem for the linearized equation (1.23) with $L(v)$ defined by (1.22) (1.14) (1.16) (1.17) (1.20) for a given $v$, with initial time $t_0 \geq 0$. We first give a preliminary result with $t_0 > 0$ where we do not study the behaviour of the solution as $t$ tends to zero.\\

\noi {\bf Proposition 3.1.} {\it Let $\rho > 1$, let $I = (0, T]$, let $v_0 \in H^\rho$ and let $v \in L^\infty (I, H^\rho)$. Let $0 \leq \rho ' < 3/2$, let $0 < t_0 \leq T$ and let $v'_0 \in H^{\rho '}$. Then the equation (1.23) has a unique solution $v' \in {\cal C}^L(I, H^{\rho '})$ with $v'(t_0) = v'_0$. The solution satisfies 
$$\parallel v'(t) \parallel_2\ = \ \parallel v'_0 \parallel_2$$

\noi for all $t \in  I$ and is unique in ${\cal C}(I,L^2)$.}\\

\noi{\bf Sketch of proof.} The result follows from the fact that the operator $L(v)$ is bounded in $\dot{H}^\sigma$ for $0 \leq \sigma < 3/2$ for all $t \in I$ and is self-adjoint in $L^2$. In fact for any $v' \in \dot{H}^\sigma$

$$\parallel \omega^{\sigma} L(v) v'\parallel_2 \ \leq \ \left ( \parallel e^{i\varphi}; M^\sigma \parallel^2 \ \parallel B(u_c);M^\sigma \parallel \ + \ 
\parallel B(v_0);M^\sigma \parallel \right )  \parallel \omega^\sigma v'\parallel_2$$

\noi and the norms in the right hand side are estimated by (\ref{2.21e}) (\ref{2.5e}) and (\ref{2.8e}). \par\nobreak
\hfill $\sq$

We next study the boundedness and continuity properties near $t = 0$ of the generic solutions of (1.23) obtained in Proposition 3.1. In view of later applications with $\rho ' = \rho$, we henceforth restrict our attention to the case $\rho ' > 1$. \\

\noi {\bf Proposition 3.2.} {\it Let $\rho > 1$, let $I = (0, T]$ and let $v \in L^\infty(I, H^\rho ) \cap {\cal C}([0,T], L^2)$ be such that $u_c$ defined by (1.16) (1.20) with $v_0 = v(0)$ satisfy the equation (\ref{2.22e}) for some real $V \in L_{loc}^\infty (I, L^\infty )$. Let $1 < \rho ' < 3/2$ and let $v' \in {\cal C}(I, H^{\rho '})$ be a solution of the equation (1.23) in $I$. Then\par

(1) $v' \in ({\cal C}^L \cap L^\infty )(I, H^{\rho '}) \cap {\cal C}_w ([0, T], H^{\rho '}) \cap {\cal C} ([0, T], H^\sigma )$ for $0 \leq \sigma < \rho '$. \par

(2) Let $0 < \theta < (1/2) \wedge (\rho -1)$. Then for all $t \in [0, T]$, $t_1 \in I$, the following estimate holds
\beq
\label{3.1e}
\parallel \omega^{\rho '} v' (t)\parallel_2 \ \leq \  \parallel \omega^{\rho '} v'(t_1) \parallel_2 \ E(|t - t_1|)
\eeq

\noi with 
\bea
\label{3.2e}
&&E(t) = E(t,a) = \exp \left ( C(\theta ) t^\theta \ a^2 \left ( 1 + a^2 (1 + |\ell n\ t|) \right )^8 \right ) \  , \\
&&a= \ \parallel v; L^\infty (I, H^\rho ) \parallel \ .
\label{3.3e}
\eea

(3) For all $t, t_1 \in [0, T]$, the following estimate holds}
\beq
\label{3.4e}
\parallel v'(t) - v'(t_1)\parallel_2 \ \leq \ C|t-t_1|^{\rho '/2} \ a^2\left ( 1 + a^2 (1 + |\ell n| t - t_1||)\right )^6  \parallel \omega^{\rho '} v' (t \vee t_1)\parallel_2 \ .\eeq

\noi {\bf Proof.} We know already that the $L^2$ norm of $v'$ is conserved. The bulk of the proof consists in deriving the estimates (\ref{3.1e}) and (\ref{3.4e}) for $t, t_1 \in I$. We begin with (\ref{3.1e}). From (1.23) we obtain 
\beq
\label{3.5e}
t\partial_t \parallel \omega^{\rho '} v'\parallel_2^2 \ = 2 \ {\rm Im} \ <\omega^{\rho '} v', \omega^{\rho '} \ L(v) v'>\ = \ J_0 + J_1 + J_2 
\eeq

\noi where
\bea
\label{3.6e}
&&J_0 = - 2\ {\rm Im} \ < \omega^{\rho '} v', \omega^{\rho '} \ U_\varphi^* \left ( B(u_c) - B_0\right ) U_\varphi v'> \  , \\
\label{3.7e}
&&J_1 = - 2\ {\rm Im} \ < v', \left [ \omega^{2\rho '}, e^{i\varphi}\right ] \left ( U^*B_0U - B_0 \right )  e^{- i\varphi} v'> \ ,\\
&&J_2 = i \ < v', e^{i\varphi} \left ( U^*\left [ \omega^{2\rho '}, B_0 \right ] U -  \left [ \omega^{2\rho '}, B_0\right ]  \right ) e^{-i\varphi} \ v'>
\label{3.8e}
\eea

\noi with $B_0 = B(v_0)$. We estimate $J_0$ by
\bea
\label{3.9e}
|J_0| &\leq& \parallel e^{i\varphi}; M^{\rho '} \parallel^2 \ \parallel B(u_c) - B_0 ; M^{\rho '} \parallel \  \parallel \omega^{\rho '}  v'\parallel_2^2 \nn \\
&\leq& C(\theta ) t^\theta \left ( 1 + \ \parallel \nabla v_0 \parallel_2^2 \ (1 + |\ell n \ t|)\right )^8 \nn \\
&&\times \prod_\pm \parallel v; L^\infty \left ( (0, t] , \dot{H}^{1 + \theta \pm 0} \right ) \parallel \ \parallel \omega^{\rho '}  v'\parallel_2^2
\eea

\noi by (\ref{2.5e}), (\ref{2.8e}), (\ref{2.21e}) and Lemma 2.10 with $2\sigma - 1/2 - 2\theta = 3/2 \pm 0$ or equivalently $\sigma = 1 + \theta \pm 0$ which allows for $\sigma < 3/2$ for $\theta < 1/2$. In order to estimate $J_1$ and $J_2$ we use the identity
\beq
\label{3.10e}
U^*AU - A = (U^* - 1) A U + A(U - 1)
\eeq

\noi and the estimate
\beq
\label{3.11e}
\parallel (U - 1)f\parallel_2 \ \leq \ t^\mu \parallel \omega^{2\mu} f \parallel_2 
\eeq

\noi with $0 \leq \mu \leq 1$. We also use Lemma 2.4 in the case
\beq
\label{3.12e}
<v_1, [\omega^{2\rho '}, f]v_2>\ \leq \ C\parallel \omega^{3/2 + 2\mu} f \parallel_2\  \parallel \omega^{\rho ' - 2\mu} v_1 \parallel_2\ \parallel \omega^{\rho '} v_2 \parallel_2
\eeq

\noi with $0 < \mu < 1/2 (<\rho '/2)$. We obtain
\bea
\label{3.13e}
&&|J_1| \ \leq \ C\ t^{\mu} \parallel \omega^{3/2+ 2\mu}\ e^{i \varphi}  \parallel_2 \ \parallel B_0 ; M^{\rho '} \parallel \  \parallel e^{i\varphi} ; M^{\rho '} \parallel \  \parallel \omega^{\rho '} v'\parallel_2^2\nn \\
&&\leq \ C\ t^{\mu} |\ell n\ t| \ \parallel \omega^{1 + \mu}\ v_0\parallel_2^2 \ \parallel \omega^{1\pm 0}\ v_0 \parallel_2^2 \left ( 1 \ + \ \parallel \nabla v_0 \parallel_2^2 |\ell n \ t|\right )^4  \ \parallel \omega^{\rho '} v'\parallel_2^2 \nn \\
\eea

\noi by (\ref{2.5e}), (\ref{2.8e}), (\ref{2.20e}) with $\sigma = 3/2 + 2\mu$ and (\ref{2.21e}), for $0 < \mu < 1/2$. Similarly, we estimate
\bea
\label{3.14e}
&&|J_2| \ \leq \ C\ t^{\mu} \parallel \omega^{3/2+ 2\mu}\ B_0  \parallel_2 \ \parallel  e^{i\varphi} ; M^{\rho '} \parallel^2 \  \parallel \omega^{\rho '}  v'\parallel_2^2\nn \\
&&\leq \ C\ t^{\mu} \parallel \omega^{1 + \mu}\ v_0\parallel_2^2 \ \left ( 1 \ + \parallel \nabla v_0\parallel_2^2 \  |\ell n\ t|\right )^4  \  \parallel \omega^{\rho '} v'\parallel_2^2 \ .
\eea

\noi Collecting (\ref{3.9e}) (\ref{3.13e})  (\ref{3.14e}) and using the fact that all the norms of $v$ occurring therein are bounded by $a$ for $0 < \mu$, $\theta < 1/2 \wedge (\rho -1)$, we obtain
\bea
\label{3.15e}
&&t\left | \partial_t  \parallel \omega^{\rho '} v'\parallel_2^2\right | \ \leq \ \Big \{ C_1(\theta ) t^\theta \ a^2 \left ( 1 + a^2 (1 + |\ell n \ t|)\right )^8\nn \\
&&+ \ C\ t^{\mu} \ a^2 \left ( 1 + a^2 |\ell n\ t|\right )^5 \Big \} \ \parallel \omega^{\rho '} v'\parallel_2^2\ .
\eea

\noi Taking $\mu = \theta$ and integrating over time yields (\ref{3.1e}) for $t, t_1 \in I$ by an elementary computation using the fact that
$$\int_{t_1}^t dt\ t^{-1 + \theta}\ |\ell n\ t|^p \leq \int_0^{t-t_1} dt\ t^{-1 + \theta }\ |\ell n\ t|^p$$

\noi for $p \geq 1$ and $0 \leq t_1 \leq t \leq 1$.\par

We next derive the estimates (\ref{3.4e}) for $t, t_1 \in I$. From (1.23) we obtain
\beq
\label{3.16e}
\partial_t \parallel v'(t) - v'_1 \parallel_2^2 \ = 2 t^{-1}\ {\rm Im} \ <v'(t) - v'_1, L(v) v'_1>
\eeq

\noi where $v'_1 = v'(t_1)$ so that
\beq
\label{3.17e}
\left | \partial_t \parallel v'(t) - v'_1 \parallel_2\right |  \ \leq \  t^{-1}\parallel L(v) v'_1 \parallel_2 \ \leq t^{-1}(K_0 + K_1)
\eeq

\noi where
\bea
\label{3.18e}
&&K_0 = \ \parallel (B(u_c) - B_0 ) U_\varphi \ v'_1 \parallel_2 \  , \\
&&K_1 = \ \parallel (U^* B_0  U - B_0) \ e^{-i \varphi} v'_1 \parallel_2 \  .
\label{3.19e}
\eea

\noi We estimate
\bea
\label{3.20e}
K_0 &\leq&\parallel \omega^{3/2 - \rho '} ( B(u_c) - B_0 )  \parallel_2 \ \parallel  e^{i\varphi} ; M^{\rho '} \parallel \  \parallel \omega^{\rho '} v'_1\parallel_2\nn \\
&\leq& \ C(\sigma )\ t^{\theta}  \left ( 1\ + \ \parallel \nabla v_0\parallel_2^2 \  (1 + |\ell n\ t|)\right )^6  \  \parallel v; L^\infty ((0,t], \dot{H}^\sigma ) \parallel^2\ \parallel \omega^{\rho '} v'_1\parallel_2 \nn \\
\eea

\noi by  (\ref{2.8e}) (\ref{2.21e}) and Lemma 2.10 with 
$$1 < \sigma < 3/2 \quad , \qquad 0 < \theta = \rho '/2 + \sigma - 1 \leq 1 \ .$$

\noi We next estimate
\bea
\label{3.21e}
K_1 &\leq& t^{\rho '/2} \parallel B_0 ; M^{\rho '}  \parallel \ \parallel  e^{i\varphi} ; M^{\rho '} \parallel \  \parallel \omega^{\rho '} v'_1\parallel_2\nn \\
&\leq& \ C\ t^{\rho '/2}  \parallel \omega^{1\pm 0} \ v_0 \parallel_2^2 \left ( 1\ + \ \parallel \nabla v_0\parallel_2^2 \  |\ell n\ t|\right )^2  \  \parallel \omega^{\rho '} v'_1 \parallel_2
\eea

\noi by (\ref{3.10e}) (\ref{3.11e})  (\ref{2.5e}) (\ref{2.8e}) (\ref{2.21e}). Substituting (\ref{3.20e})  (\ref{3.21e}) into  (\ref{3.17e}) yields
\beq
\label{3.22e}
\left | \partial_t \parallel v'(t) - v'_1 \parallel_2 \right |\  \leq \  C\ t^{-1+\rho '/2} \ a^2 \left ( 1 + a^2(1 + |\ell n \ t|)\right )^6 \parallel  \omega^{\rho '} v'_1 \parallel_2 
\eeq

\noi from which (\ref{3.4e}) follows by integration over time for $t, t_1 \in I$. We next exploit (\ref{3.1e}) and (\ref{3.4e}) in $I$ to complete the proof of the proposition. From (\ref{3.1e}) it follows that $v' \in L^\infty (I, H^{\rho '})$. From (\ref{3.4e}) and (\ref{3.1e}) it then follows that $v'$ has a limit $v'(0)$ in $L^2$ and that (\ref{3.4e}) holds for $t, t_1 \in [0, T]$. It then follows by a standard abstract argument that $v'(0) \in H^{\rho '}$, that $v' \in {\cal C}_w ([0, T], H^{\rho '}) \cap  {\cal C}  ([0, T], H^{\sigma})$ for $0 \leq \sigma < \rho '$ and that (\ref{3.1e}) holds for all $t \in [0, T]$, $t_1 \in I$. \par\nobreak
\hfill $\sq$

We have not proved so far that $v'\in {\cal C} ([0, T], H^{\rho '})$. This is true but requires a separate argument.\\

\noi {\bf Proposition 3.3.} {\it Under the assumptions of Proposition 3.2, $v' \in {\cal C}([0, T], H^{\rho '})$ and (\ref{3.1e}) holds for all $t, t_1 \ in [0, T]$.}\\

\noi {\bf Proof.} Let $t_1 \in I$, $t_1 > 0$ and let $A(t)$ be the linear map $v'_1 \to v'(t)$ defined by Propositions 3.1 and 3.2 for $0 \leq t \leq T$, with $v'_1 = v'(t_1)$. Let $0 < \varepsilon < 3/2 - \rho '$. It follows from (\ref{3.1e}) that $A(t)$ is bounded in $H^\sigma$ for $1 < \sigma < 3/2$ and in particular satisfies the estimates
\bea
\label{3.23e}
&&\parallel A(t) v'_1 ; H^{\rho '} \parallel \ \leq \  \overline{E} \parallel  v'_1 ; H^{\rho '} \parallel \\
&&\parallel A(t) v'_1 ; H^{\rho '+ \varepsilon} \parallel \ \leq \  \overline{E} \parallel  v'_1 ; H^{\rho '+ \varepsilon} \parallel 
\label{3.24e}
\eea

\noi for all $t \in [0, T]$, for some constant $\overline{E}$ independent of $t$.\par

We decompose a function $f$ into its high and low frequency parts $f_>$ and $f_<$ according to 
\beq
\label{3.25e}
\widehat{f}_{>\atop <}(\xi ) = \chi \left ( |\xi | \ \mathrel{\mathop >_{<}}\ N \right ) \widehat{f}(\xi ) 
\eeq

\noi for some (large) $N > 0$. \par

Let now $v'$ satisfy the assumptions of Proposition 3.2. We first show that\break\noindent $\parallel v'(t)_> ; \dot{H}^{\rho '}\parallel$ tends to zero uniformly in $t$ when $N \to \infty$. For that purpose we estimate
\begin{eqnarray*}
&&\parallel v'(t)_> ; H^{\rho '}\parallel\ \equiv \  \parallel \left ( A(t) v'_1\right )_> ; H^{\rho '}\parallel \\
&&\leq \ \parallel \left ( A(t) v'_{1>}\right )_>  ; H^{\rho '}\parallel \ + \ \parallel \left ( A(t) v'_{1<}\right )_> ;  H^{\rho '}\parallel \ .
\end{eqnarray*}

\noi Now
$$\parallel \left ( A(t) v'_{1>}\right )_>  ; H^{\rho '}\parallel \ \leq \ \parallel A(t) v'_{1>} ;  H^{\rho '}\parallel\ \leq \ \overline{E} \parallel v'_{1>}  ; H^{\rho '}\parallel$$

\noi by (\ref{3.23e}) while
\begin{eqnarray*}
&&\parallel \left ( A(t) v'_{1<}\right )_>  ; H^{\rho '}\parallel \ \leq \ N^{-\varepsilon} \parallel \left ( A(t) v'_{1<}\right )_> ;  H^{\rho '+ \varepsilon}\parallel\\
&&\ \leq \ N^{-\varepsilon}\parallel A(t) v'_{1<} ;  H^{\rho '+ \varepsilon}\parallel \ \leq \ \overline{E} \ N^{-\varepsilon} \parallel v'_{1<}  ; H^{\rho '+ \varepsilon}\parallel
\end{eqnarray*}

\noi by (\ref{3.24e}), so that uniformly for $t \in [0, T]$
\beq
\label{3.26e}
\parallel v'(t)_> ; H^{\rho '}\parallel\ \leq \ \overline{E} \left ( \parallel v'_{1>} ; H^{\rho '}\parallel \ + \ N^{-\varepsilon}  \parallel v'_{1<}  ; H^{\rho '+ \varepsilon}\parallel \right ) \equiv \varepsilon_N \ .
\eeq

\noi The two terms in the second member of (\ref{3.26e}) tend to zero when $N \to \infty$ for fixed $v'_1$ by definition and by the Lebesgue dominated convergence theorem respectively.We now estimate
\begin{eqnarray*}
&&\parallel v'(t) - v'(0); H^{\rho '} \parallel \ \leq \  \parallel  (v'(t) - v'(0))_> ; H^{\rho '} \parallel \ + \  \parallel  (v'(t) - v'(0))_< ; H^{\rho '} \parallel \\
&&\leq \ 2 \varepsilon_N + (1+N)^\varepsilon \parallel v'(t) - v'(0); H^{\rho '- \varepsilon} \parallel
\end{eqnarray*}

\noi The first term in the right hand side tends to zero when $N \to \infty$ uniformly in $t$, while the second term tends to zero for fixed $N$ when $t \to 0$ since $v' \in {\cal C}([0, T], H^{\rho ' - \varepsilon})$ by Proposition 3.2.\par\nobreak 
\hfill $\sq$

We can now state the main result on the Cauchy problem for the linearized equation (1.23). \\

\noi {\bf Proposition 3.4.} {\it Let $\rho  > 1$, let $I = (0, T]$ and let $v \in L^\infty(I, H^\rho ) \cap {\cal C}([0,T], L^2)$ be such that $u_c$ defined by (1.16) (1.20) with $v_0 = v(0)$ satisfy the equation (\ref{2.22e}) for some real $V \in L_{loc}^\infty (I, L^\infty )$. Let $1 < \rho ' < 3/2$ and let $v'_0 \in H^{\rho '}$. Let $t_0 \in [0, T]$. Then there exists a unique solution $v' \in  {\cal C} ([0, T], H^{\rho '})$ of the equation (1.23) with $v'(t_0) = v'_0$. Furthermore $v'$ satisfies the estimates (\ref{3.1e}) and (\ref{3.4e}) for all $t,t_1 \in [0, T]$. The solution is actually unique in ${\cal C} ([0, T], L^2)$.}\\

\noi {\bf Proof.} For $t_0 > 0$, the result follows from Propositions 3.1, 3.2 and 3.3. For $t_0 =0$, it will be proved by a limiting procedure on $t_0$. For any $t_1 \in I$, let $v'_{t_1}$ be the solution of (1.23) with $v'_{t_1}(t_1) = v'_0$ given by Propositions 3.1 and 3.2. Let now $0 < t_1 <t_2 \leq T$. It follows from (\ref{3.1e}) that 
\beq
\label{3.27e}
\parallel \omega^{\rho '} v'_{t_i} (t)\parallel_2 \ \leq \ E \left ( |t - t_i| \right ) \parallel \omega^{\rho '} v'_0\parallel_2 
\eeq

\noi for $i = 1,2$ and for all $t \in [0, T]$. Furthermore, from  (\ref{3.4e}) and  (\ref{3.27e}) and from $L^2$-norm conservation, it follows that
\begin{eqnarray}
\label{3.28e}
&&\parallel v'_{t_2}(t) - v'_{t_1}(t) \parallel_2  \ = \  \parallel  v'_{t_2}(t_1) - v'_0\parallel_2 \ = \  \parallel  v'_{t_2}(t_1)  - v'_{t_2}(t_2) \parallel_2\nn  \\
&&\leq \ C|t_2 - t_1|^{\rho '/2} \ a^2 \left ( 1 + a^2 (1 + |\ell n (t_2 - t_1)|)\right )^6 \ E(t_2)   \parallel \omega^{\rho '} v'_0 \parallel_2 \ .
\end{eqnarray}

\noi From  (\ref{3.28e}) it follows that $v'_{t_1}$ converges in $L^\infty (I, L^2)$ norm to some $v' \in {\cal C}([0,T], L^2)$ when $t_1 \to 0$. From the uniform estimate (\ref{3.27e}) it follows by abstract arguments that $v' \in ({\cal C}_w \cap L^\infty )([0,T], H^{\rho '}) \cap {\cal C} ([0, T], H^\sigma )$ for $0 \leq \sigma < \rho '$, that $v'$ satisfies the estimates of Proposition 3.2 and that $v'(0) = v'_0$. Furthermore $v'$ is easily seen to satisfy (1.23) in $I$, so that $v' \in {\cal C}^L(I, H^{\rho '})$. It remains to be proved that actually $v'$ is strongly continuous in $H^{\rho '}$ at $t = 0$. This follows from Proposition 3.3, which has not been used so far. Alternatively it follows from the estimate (\ref{3.27e}) with $t_1 = 0$ that 
$$\lim_{t\to 0} \ \sup \parallel \omega^{\rho '}v'(t)\parallel_2\ \leq \ \parallel \omega^{\rho '}v'_0\parallel_2\ E(0) = \ \parallel \omega^{\rho '}v'(0)\parallel_2$$

\noi which together with weak continuity implies strong continuity at $t=0$. \par\nobreak \hfill $\sq$

\noi {\bf Remark 3.1.} Note that in the case where $t_0 = 0$, Proposition 3.3 is not needed for the proof of Proposition 3.4.

\mysection{The nonlinear Cauchy problem at time zero\break\noindent for v}
\hspace*{\parindent} 
In this section we prove that the non linear equation (1.21) for $v$ with initial data at time $t_0$ has  a unique solution in a small time interval by showing that the map $\Gamma$~; $v \to v'$ defined by Proposition 3.4 with $t_0 = 0$ is a contraction. For that purpose, we need to estimate the difference of two solutions of the linearized equation (1.23). \\

\noi {\bf Lemma 4.1.} {\it Let $\rho > 1$, let $0 < \theta < (1/2) \wedge (\rho - 1)$, let $I = (0,T]$ and let $v_i$, $i = 1,2$ satisfy the assumptions of Proposition 3.4 with $v_i(0) = v_0 \in H^\rho$. Let $1 < \rho ' < 3/2$ and let $v'_i$, $i =1,2$ be the solutions of the equation (1.23) with $v'_i (0) = v'_0 \in H^{\rho '}$ obtained in Proposition 3.4. Let $v_- = v_2 - v_1$ and $v'_- = v'_2 - v'_1$. Then the following estimate holds for all $t$, $0 < t \leq T$~:
\beq
\label{4.1e}
\parallel v'_-; L^\infty ((0, t], H^{\rho '}) \parallel \ \leq \ C\ t^\theta E(t,a) aa' \left ( 1 + a^2 (1 + |\ell n\ t|)\right )^8  \parallel v_-; L^\infty ((0, t], H^\rho ) \parallel
\eeq 

\noi where $E(t, a)$ is defined by (\ref{3.1e}) and} 
$$a = \ {\rm Max}\parallel v_i; L^\infty ((0, T], H^\rho ) \parallel\ , \quad a' = \ {\rm Max}\parallel v'_i; L^\infty ((0, T], H^{\rho '} ) \parallel\  .$$

\noi {\bf Proof.} From (1.23) we obtain
\begin{eqnarray*}
i\partial_t v'_- &=&t^{-1} \left ( L_2 v'_2 - L_1 v'_1 \right ) \\
&=&t^{-1} \left ( L_2 v'_- + L_- v'_1 \right ) 
\end{eqnarray*}

\noi where $L_i =L(v_i)$ and
\begin{eqnarray*}
&&L_- = L_2 - L_1 = - U_\varphi^* B_- U_\varphi \ , \\
&&B_- = B\left ( u_{c2} \right ) - B\left ( u_{c1} \right ) \ .
\end{eqnarray*}

\noi We estimate for $0 \leq \sigma ' < 3/2$
\beq
\label{4.2e}
t\partial_t\ \parallel \omega^{\sigma '} v'_- (t)\parallel_2^2 \ = \ 2\ {\rm Im}  \left (  < \omega^{\sigma '}  v'_-, \omega^{\sigma '} L_2 v'_- > \ + \ <   \omega^{\sigma '} v'_-, \omega^{\sigma '} L_- v'_1 > \right )\ .
\eeq

\noi By the estimates in the proof of Proposition 3.2 (see in particular  (\ref{3.15e})) we obtain
\beq
\label{4.3e}
\parallel \omega^{\sigma '} v'_- (t)\parallel_2 \ \leq \ E(t,a) \int_0^t dt'\ t{'}^{-1} \parallel  \omega^{\sigma '} L_- v'_1(t')\parallel_2 \ .
\eeq

\noi We next estimate
\bea
\label{4.4e}
&&\parallel \omega^{\sigma '} L_- v'_1 \parallel_2 \ \leq \ \parallel e^{i \varphi}; M^{\sigma '} \parallel^2 \ \parallel B_-;  M^{\sigma '}\parallel \ \parallel \omega^{\sigma '} v'_1\parallel_2 \nn \\
&&\leq C (1 + a^2|\ell n\ t|)^4 \ \parallel \omega^{3/2 \pm 0} B_-  \parallel_2 \ \parallel  \omega^{\sigma '} v'_1\parallel_2 \ .
\eea 

\noi by (\ref{2.5e}) (\ref{2.21e}). From (1.14) we obtain
\beq
\label{4.5e}
B_-(t) = \int_1^\infty d\nu \ \nu^{-3} \omega^{-1} \sin (\omega (\nu - 1)) D_0 (\nu ) (|u_{c2}|^2 - |u_{c1}|^2) (t/ \nu ) \ .
\eeq

\noi By the same method as in the proof of (\ref{2.8e}), we estimate 
\beq
\label{4.6e}
\parallel \omega^{2\sigma - 1/2} B_-(t) \parallel_2 \ \leq \ C \int_1^\infty d\nu \ \nu^{-2\sigma} \parallel \omega^{\sigma} u_{c+} (t/\nu ) \parallel_2 \ \parallel \omega^{\sigma} u_{c-}(t/\nu ) \parallel_2                                                                                                                                                                                                                                                                                                                                                                                                                                                                                                                                                                                                                                                                                                                                                                                                                                                                                                                                                                                                                                                                                                                                                                                                                                                                                                                                                                                                                                                                                                                                                                                                                                                                                                                                                                                                                                                                                                                                                                                                                                                                                                                            \eeq

\noi for $1/2 < \sigma < 3/2$, with $u_{c_\pm} = u_{c2} \pm u_{c1}$. On the other hand, in the same way as in Lemma 2.6, the local conservation law (\ref{2.15e}) for $u_{ci}$ implies
\beq
\label{4.7e}
\partial_t\ <|u_{c2}|^2 - |u_{c1}|^2, \psi > \ = - (i/4) \ (<u_{c+}, [\Delta , \psi ] u_{c-}> \ + \ <u_{c-}, [\Delta , \psi ] u_{c+} > ) 
\eeq

\noi for any test function $\psi$, so that
\bea
\label{4.8e}
&&\parallel \omega^{2\sigma - 7/2} \left ( |u_{c2}(t)|^2 - |u_{c1}(t)|^2 - \left ( |u_{c2}(t_1)|^2 - |u_{c1}(t_1)|^2\right ) \right )  \parallel_2\nn \\
&& \ \leq \ C \int_{t_1}^t dt' \parallel  \omega^{\sigma} u_{c+}(t')\parallel_2 \ \parallel \omega^{\sigma} u_{c-}(t') \parallel_2 
\eea

\noi for $1/2 < \sigma < 3/2$. Substituting (\ref{4.8e}) into the definition of $B_-$ yields 

$$\parallel \omega^{2\sigma - 5/2} \left ( B_-(t) - B_-(t_1)\right )  \parallel_2 \ \leq  \int_1^\infty d\nu \ \nu^{2 -2\sigma} \int_{t_1/\nu}^{t/\nu} dt'\parallel \omega^{\sigma} u_{c+} (t' ) \parallel_2 \ \parallel \omega^{\sigma} u_{c-} (t' ) \parallel_2$$
\beq
\label{4.9e}
\leq \ C \ t \int_1^\infty d\nu \ \nu^{1-2\sigma} \parallel \omega^{\sigma} u_{c+} (t/\nu ) \parallel^2 \ \parallel \omega^{\sigma} u_{c-}(t/\nu ) \parallel_2                                                                                                                                                                                                                                                                                                                                                                                                                                                                                                                                                                                                                                                                                                                                                                                                                                                                                                                                                                                                                                                                                                                                                                                                                                                                                                                                                                                                                                                                                                                                                                                                                                                                                                                                                                                                                                                                                                                                                                                                                                                                                                                                                                                                                                                                                  \eeq

\noi for $1 < \sigma < 3/2$, in the same way as in the proof of Lemma 2.7. We know from Lemma 2.10 that $B_- (t_1)$ tends to $B_- (0)$ when $t_1 \to 0$ in the norms appearing in (\ref{4.9e}) and that $B_- (0) = 0$ since $v_1(0) = v_2(0)$. Therefore 
\beq
\label{4.10e}
\parallel \omega^{2\sigma - 5/2} B_-(t) \parallel_2 \ \leq \ C \ t \int_1^\infty d\nu \ \nu^{1-2\sigma} \parallel \omega^{\sigma} u_{c+} (t/\nu ) \parallel_2 \ \parallel \omega^{\sigma} u_{c-}(t/\nu ) \parallel_2   \ .                                                                                                                                                                                                                                                                                                                                                                                                                                                                                                                                                                                                                                                                                                                                                                                                                                                                                                                                                                                                                                                                                                                                                                                                                                                                                                                                                                                                                                                                                                                                                                                                                                                                                                                                                                                                                                                                                                                                                                                                                                                                                                                         \eeq

\noi Interpolating between (\ref{4.6e}) and (\ref{4.10e}) we obtain 
$$ \parallel \omega^{2\sigma - 1/2- 2\theta} B_-(t) \parallel_2 \ \leq \ C \ t^\theta \int_1^\infty d\nu \ \nu^{1-2\sigma} \parallel \omega^{\sigma} u_{c+} (t/\nu ) \parallel_2 \ \parallel \omega^{\sigma} u_{c-}(t/\nu ) \parallel_2$$                                                                                                                                                                                                                                                                                                                                                                                                                                                                                                                                                                                                                                                                                                                                                                                                                                                                                                                                                                                                                                                                                                                                                                                                                                                                                                                                                                                                                                                                                                                                                                                                                                                                                                                                                                                                                                                                                                                                                                                                                                                                                                                         

\noi for $1 < \sigma < 3/2$ and $0 \leq \theta \leq 1$, 
 \beq
\label{4.11e}
\cdots  \leq C\ t^\theta  (1 + a^2(1 + |\ell n\ t|))^4 a\ \parallel v_-; L^\infty ((0, t], \dot{H}^\sigma )  \parallel 
\eeq 

\noi (see (\ref{2.25e}) (\ref{2.26e})). Substituting (\ref{4.11e}) into (\ref{4.4e}) with $\sigma = 1 + \theta \pm 0$ and substituting the result into (\ref{4.3e}) yields (\ref{4.1e}). \par \nobreak \hfill $\sq$

We can now state the main result on the Cauchy problem at time zero for the equation (1.21).\\

\noi {\bf Proposition 4.1.} {\it Let $1 < \rho < 3/2$ and let $v_0 \in H^\rho$. Then there exists $T > 0$ and there exists a unique solution $v \in {\cal C}([0, T], H^\rho )$ of the equation (1.21) with $v(0) = v_0$. One can ensure that 
 \bea
\label{4.12e}
&& \parallel v; L^\infty ([0, T], H^\rho )  \parallel \ \leq \ R=2 \parallel v_0; H^\rho \parallel \\
&&C\ T^\theta R^2 \left ( 1 + R^2 (1 + |\ell n\ T|) \right )^8 = 1/2
\label{4.13e}
\eea

\noi for some $\theta$, $0 < \theta < \rho - 1$, and some $C$ independent of $v_0$.}\\

\noi {\bf Proof.} Let $v_0 \in H^\rho$, let $B_0 = B(v_0)$ and $\varphi = - \ell n\ t \ B_0$. Let $T > 0$. For any $v\in {\cal C}([0, T], H^\rho )$, we define $u_c = U_\varphi v$. Let $F(T, v_0)$ be the set of $v \in  {\cal C}([0, T], H^\rho )$ such that $v(0) = v_0$ and that $u_c$ satisfy the equation (\ref{2.22e}) in $(0, T]$ for some real $V \in L_{loc}^\infty ((0, T], L^\infty )$. It follows from Proposition 3.4 that $F(T, v_0)$ is stable under the map $\Gamma : v\to v'$ defined by that proposition with $t_0 = 0$ and $v'_0 = v_0$. Let $B(R)$ be the ball of radius $R$ in ${\cal C}([0, T], H^\rho$). From Proposition 3.2 it follows that $B(R) \cap F(T, v_0)$ is stable under $\Gamma$ if
\beq
\label{4.14e}
E(T, R) \leq 2
\eeq 

\noi with $R = \ 2 \parallel v_0; H^\rho \parallel$. Furthermore by Lemma 4.1, $\Gamma$ is a contraction in the $L^\infty ([0, T]), H^\rho )$ norm on that set if
\beq
\label{4.15e}
C\ T^\theta E(T,R)\ R^2 \left ( 1 + R^2(1 + |\ell n \ T|) \right )^8 \leq 1/2\ .
\eeq

\noi Therefore for $T$ sufficiently small to satisfy (\ref{4.14e}) (\ref{4.15e}), namely under a condition of the type (\ref{4.13e}), the map $\Gamma$ has a unique fixed point in $B(R)$ provided $F(T, v_0)$ is non empty. The set $F(T,v_0)$ is non empty because it contains the solution of the linear equation 
$$i \partial_t v^{(0)} = t^{-1} \left ( U_\varphi^* \ B_0 \ U_\varphi - B_0 \right ) v^{(0)}$$

\noi with $v^{(0)}(0) = v_0$ obtained by a simplified version of Proposition 3.2. Clearly the fixed point $v$ satisfies the equation (1.21) and therefore belongs to $F(T, v_0)$. \par\nobreak \hfill $\sq$

\mysection{Continuity with respect to initial data}
\hspace*{\parindent} 
In this section we prove that the map $v_0 \to v$ defined in Proposition 4.1 is continuous in the natural norms. For that purpose, we need to estimate the difference of two solutions of the linearized equation (1.23) corresponding to two functions $v_1$ and $v_2$ not necessarily satisfying the condition $v_1(0) = v_2 (0)$. The following lemma is an extension of Lemma 4.1 where we drop that assumption. \\

\noi {\bf Lemma 5.1.} {\it Let $\rho > 1$, let $I = (0, T]$ and let $v_i$, $i = 1, 2$ satisfy the assumptions of Proposition 3.4 with $v_i (0) = v_{0i} \in H^\rho$. Let $1 <\rho ' < 3/2$ and let $v'_i$, $i =1,2$, be the solutions of the equation (1.23) with $v'_i (0) = v'_{0i}  \in H^{\rho '}$ obtained in Proposition 3.4. Let $v_- = v_2 - v_1$, $v_{0-} = v_{02} - v_{01}$ and similarly $v'_- = v'_2 - v'_1$, $v'_{0-} = v'_{02}  - v'_{01}$. Let $0 < 2\mu \leq \rho '$ and $\theta = \mu + (\rho - 1)/2 \wedge (1/8)$. Then the following estimate holds for all $t$, $0 < t \leq T$~: 
\bea
\label{5.1e}
&&\parallel \omega^{\rho ' - 2  \mu} \ v'_-(t) \parallel_2\ \leq E(t,a) \Big ( \parallel \omega^{\rho ' - 2  \mu} \ v'_{0-} \parallel_2\nn \\
&&+ \ C(\rho ) aa' \left ( \mu^{-7} t^\mu L^{10}  \parallel  v_{0-} ; H^\rho \parallel\ + \ t^\theta L^8  \parallel v_- ; L^\infty ((0, t], H^\rho )  \parallel \right ) \Big )  
\eea

\noi where $E(t,a)$ is defined by (\ref{3.1e}), 
\begin{eqnarray*}
&&a = \ {\rm Max} \parallel v_i ; L^\infty ((0, T], H^\rho ) \parallel\ , a' =  \ {\rm Max} \parallel v'_i ; L^\infty ((0, T], H^{\rho '}) \parallel\ , \\
&&L = 1 + a^2 (1 + |\ell n\ t|) \  .
\end{eqnarray*}

\noi The constant $C(\rho )$ depends on $\rho$ but is independent of $\mu$.}\\

\noi {\bf Proof.} In the same way as in the proof of Lemma 4.1 we obtain from (1.23)
\begin{eqnarray*}
i \partial_t v'_-&=&t^{-1} \left ( L_2 v'_2 - L_1 v'_1 \right ) \\
&=&t^{-1} \left ( L_2 v'_-  + L_- v'_1 \right ) 
\end{eqnarray*}

\noi where $L_i = L(v_i)$ and $L_- = L_2 - L_1$, from which we obtain (\ref{4.2e}) and (\ref{4.3e}). Now however $L_-$ is more complicated because we do not assume that $v_1(0) = v_2(0)$. Let $B_{0i} = B(v_i (0))$, $B_{0-} = B_{02} - B_{01}$ and correspondingly $\varphi_i = - (\ell n\ t)B_{0i}$, $\varphi_- = \varphi_2 - \varphi_1$. Let $B_i = B(u_{ci})$ and $B_- = B_2 - B_1$. Then 
\bea
\label{5.2e}
-L_- &=&U_{\varphi_2}^* B_2 U_{\varphi_2} - B_{02} - U_{\varphi_1}^* B_1 U_{\varphi_1} + B_{01}\nn \\
&=&U_{\varphi_2}^* (B_2 - B_{02}) U_{\varphi_2} - U_{\varphi_1}^* (B_1 - B_{01}) U_{\varphi_1} \nn \\
&&+ e^{i\varphi_2} \left ( U^* B_{02} U - B_{02}\right ) e^{-i \varphi_2} - e^{i \varphi_1} \left ( U^* B_{01} U - B_{01}\right ) e^{-i \varphi_1} \nn \\
&=& J_0 + J_1 + J_2 + J_3
\eea

\noi where 
\bea
\label{5.3e}
J_0 &=& \left ( e^{i\varphi_2} - e^{i\varphi_1}\right ) U^* (B_2 - B_{02}) U_{\varphi_2} + U_{\varphi_1}^* (B_2 - B_{02}) \left ( e^{-i\varphi_2} - e^{-i\varphi_1}\right )\ , \\
\label{5.4e}
J_1 &=& U_{\varphi_1}^* (B_-  - B_{0-})  U_{\varphi_1} \ , \\
\label{5.5e}
J_2 &=&\left ( e^{i\varphi_2} - e^{i\varphi_1}\right ) \left ( U^* B_{02} U- B_{02} \right ) e^{-i\varphi_2} \nn \\
&&+ e^{i\varphi_1} \left ( U^* B_{02} U - B_{02} \right ) \left ( e^{-i\varphi_2} - e^{-i\varphi_1}\right ) \ , \\
J_3 &=& e^{i\varphi_1} \left ( U^* B_{0-} U - B_{0-} \right  ) e^{-i\varphi_1} \  .
\label{5.6e}
\eea 

\noi We now take $0 \leq \sigma ' = \rho ' - 2 \mu < \rho ' (<3/2)$ and we estimate the previous quantities as follows
\bea
\label{5.7e}
&&\parallel \omega^{\sigma '} J_0 v'_1 \parallel_2 \ \leq \ C \left ( \prod_i \parallel e^{i\varphi_i}; M^{\rho '} \parallel \right ) \parallel \left ( e^{i\varphi_-} - 1 \right ) ; M^{\rho '} \parallel \nn \\
&&\times \prod_\pm \parallel \omega^{3/2 - 2 \mu _\pm} (B_2 - B_{02})  \parallel_2^{1/2} \   \parallel \omega^{\rho '} v'_1 \parallel_2 \ ,                                                                                                                                                                                                                                                                                                                                                                                                                                                                                                                                                                                                                                                                                                                                                                                                                                                                                                                                                                                                                                                                                                                                                                                                                                                                                                                                                                                                                                                                                                                                                                                                                                                                                                                                                                                                                                                                                                                                                                                                                                                                                                                         
\eea
\beq
\parallel \omega^{\sigma '} J_1 v'_1 \parallel_2 \ \leq \ C\parallel e^{i\varphi_1}; M^{\rho '} \parallel^2 \prod_\pm \parallel \omega^{3/2 - 2 \mu _\pm} (B_- - B_{0-})  \parallel_2^{1/2} \   \parallel \omega^{\rho '} v'_1 \parallel_2 \ ,  \nn \\                                                                                                                                                                                                                                                                                                                                                                                                                                                                                                                                                                                                                                                                                                                                                                                                                                                                                                                                                                                                                                                                                                                                                                                                                                                                                                                                                                                                                                                                                                                                                                                                                                                                                                                                                                                                                                                                                                                                                                                                                                                                                                                       
\label{5.8e} 
\eeq

\noi with
$$\left \{ \begin{array}{ll} \mu_\pm = \mu &\hbox{for $\mu \geq \varepsilon$} \ ,\\ \\ \mu_\pm = \mu \pm \varepsilon &\hbox{for $\mu < \varepsilon$} \end{array}\right .$$

\noi for some fixed $\varepsilon$ and with constants $C$ which depend on $\varepsilon$ but can be taken independent of $\mu$. Here we have used estimates of the type
\beq
\label{5.9e}
\parallel \omega^{\sigma '} Bv \parallel_2 \ \leq \ C \prod_\pm \parallel \omega^{3/2-2\mu_\pm} B \parallel_2^{1/2} \ \parallel \omega^{\rho '} v \parallel_2
\eeq 

\noi which follow from a minor variation of Lemma 2.3, with constants $C$ satisfying the property quoted above. \par

We next estimate
\beq
\label{5.10e}
\parallel \omega^{\sigma '} J_2 v'_1 \parallel_2 \ \leq \ C\ t^{\mu} \left ( \prod_i \parallel e^{i\varphi_i}; M^{\rho '} \parallel \right )  \parallel \left ( e^{i\varphi_-}- 1 \right ) ; M^{\rho '} \parallel \parallel B_{02}; M^{\rho '}\parallel \   \parallel \omega^{\rho '} v'_1 \parallel_2  \  , 
\eeq
\beq
\parallel \omega^{\sigma '} J_3 v'_1 \parallel_2 \ \leq \ C\ t^{\mu} \parallel e^{i\varphi_1}; M^{\rho '} \parallel^2 \ \parallel B_{0-}; M^{\rho '}  \parallel \   \parallel \omega^{\rho '} v'_1 \parallel_2    
\label{5.11e}
\eeq 

\noi where we have used (\ref{3.11e}) and where the constants are again independent of $\mu$. Collecting (\ref{5.7e}) (\ref{5.8e}) (\ref{5.10e}) (\ref{5.11e}), we obtain
\beq 
\label{5.12e}
\parallel \omega^{\sigma '} L_- v'_1 \parallel_2 \ \leq \ C\ \left ( \mathrel{\mathop {\rm Max}_i} \parallel e^{i\varphi_i}; M^{\rho '} \parallel^2 \right ) K \parallel  \omega^{\rho '} v'_1 \parallel_2    
\eeq 

\noi where
\bea
\label{5.13e}
&&K = \ \parallel ( e^{i\varphi_-}- 1); M^{\rho '} \parallel \left ( \prod_\pm \parallel  \omega^{3/2 - 2 \mu _\pm} (B_2 - B_{02} \parallel_2^{1/2} \ + \  t^\mu  \parallel B_{02} ; M^{\rho '} \parallel  \right )\nn \\
&&+ \  \prod_\pm \parallel \omega^{3/2 - 2 \mu _\pm} (B_- - B_{0-})  \parallel_2^{1/2} \   + t^\mu \parallel B_{0-} ; M^{\rho '} \parallel \ .                                                                                                                                                                                                                                                                                                                                                                                                                                                                                                                                                                                                                                                                                                                                                                                                                                                                                                                                                                                                                                                                                                                                                                                                                                                                                                                                                                                                                                                                                                                                                                                                                                                                                                                                                                                                                                                                                                                                                                                                                                                                                                                       
\eea

\noi We have already estimated
\bea
\label{5.14e}
&&\parallel e^{i\varphi_i}; M^{\rho '} \parallel \ \leq \ C \left ( 1 + \parallel \nabla v_{0i} \parallel_2^2 \ |\ell n \ t|\right )^2 \ \leq \ C\ L^2 \\
&&\parallel B_{02}; M^{\rho '}\parallel \   \ \leq \ C\ a^{2} \ .    
\label{5.15e}
\eea 

\noi Furthermore, by (\ref{2.5e}) (\ref{2.8e}) (\ref{2.21e}), 
\bea
\label{5.16e}
&&\parallel B_{0-} ; M^{\rho '} \parallel \ \leq \ C\parallel \omega^{1\pm 0} v_{0+} \parallel_2 \ \parallel \omega^{1\pm 0} v_{0-} \parallel_2\ \leq \ C\ a \parallel v_{0-}; H^\rho  \parallel \ , \\
&&\parallel (e^{i\varphi_-} - 1); M^{\rho '} \parallel \ \leq \ C\Big ( \parallel \varphi_- \parallel_\infty \ + \ \parallel \omega^{3/2} \varphi_-\parallel_2 \left ( 1\ + \ 
\parallel \omega^{3/2} \varphi_-\parallel_2 \right ) \nn \\
&&\leq \ C|\ell n\ t| \Big \{  \parallel B_{0-}  \parallel_\infty \ + \  \parallel \omega^{3/2} B_{0-}  \parallel_2  \left ( 1\ + \  |\ell n \ t| \parallel \omega^{3/2}  B_{0-} \parallel_2\right ) \nn \\
&&\leq \ C|\ell n\ t| \parallel \omega^{1\pm 0} v_{0+}  \parallel_2 \  \parallel \omega^{1\pm 0} v_{0-}  \parallel_2  \left ( 1\ + \  |\ell n \ t| \parallel \nabla v_{0+}   \parallel_2 \ \parallel \nabla v_{0-} \parallel_2 \right )\ . \nn \\
&&\leq \ C\ a |\ell n\ t| (1 + a^2 |\ell n\ t |) \parallel v_{0-}; H^\rho  \parallel \ . 
\label{5.17e}
\eea

\noi We next estimate
\beq 
\label{5.18e}
\parallel  \omega^{3/2 - 2 \mu _\pm} (B_2 - B_{02} \parallel_2  \  \leq \ C\ t^{\theta_\pm}  a^2\ L^4
\eeq

\noi by Lemma 2.10, especially (\ref{2.23e}), with 
$$0 < \theta_\pm = \mu_\pm + \sigma - 1 \leq 1$$

\noi for some $\sigma$ satisfying $1 < \sigma < 3/2$, $\sigma \leq \rho$. We choose 
$$2\varepsilon = \sigma - 1 = (\rho -1) \wedge 1/4$$

\noi so that
$$\mu + (\sigma - 1)/2 \equiv \theta = \theta_- < \theta_+ \leq 1$$

\noi and
\beq 
\label{5.19e}
\parallel  \omega^{3/2 - 2 \mu _\pm} (B_2 - B_{02} \parallel_2  \  \leq \ C\ t^{\theta}  a^2\ L^4
\eeq

\noi for $0 \leq 2\mu < \rho '$, with a constant $C$ depending only on $\rho$. \par

It remains to estimate the last but one term in (\ref{5.13e}). Taking the limit $t_1 \to 0$ in (\ref{4.9e}) and using again Lemma 2.10 with $\theta = 1$, we obtain
\beq
\label{5.20e}
\parallel  \omega^{2\sigma - 5/2} (B_-(t)  - B_{0-}) \parallel_2  \  \leq \ C\ t \int_1^\infty  d\nu \ \nu^{1-2\sigma} \parallel \omega^\sigma u_{c+} (t/\nu )\parallel_2\   \parallel \omega^\sigma u_{c-} (t/\nu )\parallel_2
\eeq 

\noi which differs from (\ref{4.10e}) by the fact that now $B_{0-} \not= 0$. Interpolating between (\ref{5.20e}) and (\ref{4.6e}) together with the simpler estimate 
\beq
\label{5.21e}
\parallel  \omega^{2\sigma - 1/2} B_{0-} \parallel_2  \  \leq \ C\ a \parallel \omega^\sigma v_{0-} \parallel_2
\eeq 

\noi we obtain 
\bea
\label{5.22e}
&&\parallel  \omega^{2\sigma - 1/2-2 \theta} (B_-(t)  - B_{0-}) \parallel_2  \  \leq \ C\ t^\theta \int_1^\infty  d\nu \ \nu^{1-2\sigma}\nn \\
&&\times \left ( \parallel \omega^\sigma u_{c+} (t/\nu )\parallel_2\   \parallel \omega^\sigma u_{c-} (t/\nu )\parallel_2 \ + \ a  \parallel \omega^\sigma v_{0-} \parallel_2\right )
\eea 

\noi for $1 < \sigma < 3/2$, $\sigma \leq \rho$ and $0 < \theta \leq 1$. We estimate $ \parallel \omega^\sigma u_{c+}\parallel_2$ by (\ref{2.25e}). On the other hand, from
$$u_{c-}  = e^{-i \varphi_1} \left ( \left ( e^{-i\varphi_-} - 1\right ) v_2 + v_- \right  )$$

\noi we obtain
\beq
\label{5.23e}
\parallel  \omega^{\sigma}u_{c-}  \parallel_2  \  \leq \ C \parallel e^{-i \varphi_1} ; M^\sigma \parallel \left (  \parallel (e^{-i\varphi_-} - 1); M^\sigma \parallel \ 
\parallel \omega^\sigma v_2 \parallel_2 \ + \ \parallel \omega^\sigma v_- \parallel_2\right )\ .
\eeq 

\noi Substituting (\ref{5.14e}) (\ref{5.17e}) into (\ref{5.23e}) yields
\beq
\label{5.24e}
\parallel  \omega^{\sigma}u_{c-}  \parallel_2  \  \leq \ C L^2 \left ( a^2|\ell n\ t|L \parallel v_{0-} ; H^\rho \parallel \\ + \   \parallel v_- ; L^\infty ((0, t], H^\rho  )\parallel \right )\  . 
\eeq

\noi Substituting (\ref{5.24e}) into (\ref{5.22e}) yields (see (\ref{2.25e}))
\beq
\label{5.25e}
\parallel \omega^{3/2 - 2 \mu _\pm} (B_- - B_{0-})  \parallel_2 \   \leq \ C\  t^\theta\ a \left ( L^6 \parallel v_{0-} ; H^{\rho } \parallel\ + \  L^4 \parallel v_{-} ; L^{\infty} ((0, t], H^\rho )\parallel \right ) 
\eeq 

\noi with the same $\theta$ and with $C$ depending only on $\rho$ as in (\ref{5.19e}).\par

Substituting  (\ref{5.15e})-(\ref{5.19e}) and (\ref{5.25e}) into (\ref{5.13e}) yields
\beq
\label{5.26e}
K \leq C\ a \left ( t^\mu L^2  \parallel  v_{0-} ; H^{\rho}\parallel\ + \ t^\theta \left ( L^6  \parallel v_{0-} ; H^{\rho} \parallel\ + \  L^4 \parallel v_{-} ; L^{\infty} ((0, t], H^\rho )\parallel \right ) \right ) \ . 
\eeq

\noi Substituting (\ref{5.26e}) into (\ref{5.12e}), substituting the result into (\ref{4.3e}), integrating over time and using the fact that $\theta \geq \mu \vee ((\rho - 1)/2 \wedge (1/8))$ and that
$$\int_0^t dt'\ t{'}^{\lambda - 1} |\ell n \ t'|^p \ \leq C\ \lambda^{-(p+1)}  t^{\lambda} (1 + |\ell n\ t|)^p$$

\noi yields (\ref{5.1e}). \par \nobreak \hfill $\sq$

We can now prove the continuity of the map $v_0 \to v$ defined in Proposition 4.1.\\

\noi {\bf Proposition 5.1.} {\it Let $1 < \rho < 3/2$, let $R > 0$ and let $T$ be defined by (\ref{4.13e}) with $\theta = (\rho - 1)/2 \wedge (1/8)$. Let $B(R)$ be the ball of radius $R$ in $H^\rho$.\par

(1) Let $1 < \lambda < \rho$. Then the map $v_0 \to v$ defined by Proposition 4.1 is continuous from $H^\lambda$ to $L^\infty ((0, T], H^\lambda )$ uniformly for $v_0 \in B(R/2)$. Furthermore, there exists $T_\lambda$, $0 < T_\lambda \leq T$, such that for two solutions $v_i$, $i = 1,2$ of the equation (1.21) with $v_i (0) = v_{0i} \in B(R/2)$ as obtained in Proposition 4.1, the following estimate holds~:
\bea
\label{5.27e}
\parallel  v_{-} ; L^{\infty}((0, T_\lambda ], H^\lambda )\parallel &\leq& 2E (T_\lambda , R) \left ( 1 + C(\lambda ) R^2(\rho - \lambda )^{-7} \ T_\lambda^{(\rho - \lambda )/2} \ L (T_\lambda )^{10} \right )  \nn \\
&&\times \parallel v_{0-} ; H^{\lambda} \parallel 
\eea

\noi where $v_- = v_2 - v_1$ and $v_{0-} = v_{02} -  v_{01}$. One can define $T_\lambda$ by 
\beq
\label{5.28e}
C(\lambda )\ E(T_\lambda , R) R^2 \ T_\lambda^\theta \ L(T_\lambda )^8 = 1/2
\eeq

\noi where $E(T, \cdot )$ is defined by (\ref{3.2e}) and
$$L(t) = 1 + R^2 (1 + |\ell n\ t|) \ .$$

(2) The map $v_0 \to v$ defined by Proposition 4.1 is (pointwise) continuous from $H^\rho$ to $L^\infty ((0, T] , H^\rho )$ for $v_0 \in B(R/2)$.}\\

\noi {\bf Proof}

\noi \underline{Part (1)}. From (\ref{5.1e}) we obtain
\bea
\label{5.29e}
&&\parallel  v'_{-} ; L^{\infty}((0, t ]; H^{\rho ' - 2 \mu} )\parallel \ =  \ \mathrel{\mathop {\rm Sup}_{0 < t' \leq t}}\  \mathrel{\mathop {\rm Sup}_{\mu \leq \mu ' \leq \rho ' /2}} \parallel \omega^{\rho ' - 2 \mu '} v'_-(t') \parallel_2 \nn \\
&&\leq E(t,a) \Big ( \parallel  v'_{0-} ; H^{\rho ' - 2 \mu} \parallel \ + \ C(\rho ) aa' \Big ( \mu^{-7} t^\mu L^{10} \parallel  v_{0-} ; H^{\rho} \parallel \nn \\
&&+ \ t^\theta L^\infty \parallel  v_{-} ; L^{\infty}((0, t ], H^\rho)\parallel \Big ) \Big ) \ .
\eea

\noi Changing $(\rho ' , \rho )$ to $(\rho , \lambda )$ satisfying $1 < \lambda = \rho - 2 \mu < \rho$ and using the fact that $v'_i(t) = v_i (t) \in B(R)$ for all $t\in (0, T]$, we obtain 
\bea
\label{5.30e}
&&\parallel  v_{-} ; L^{\infty}((0, t ]; H^\lambda  )\parallel \  \leq E(t,R) \Big ( \left ( 1 + C(\lambda ) R^2 \mu^{-7} t^\mu L^{10} \right )  \parallel  v_{0-} ; H^\lambda  \parallel \nn \\
&&+ \ C(\lambda ) R^2  \ t^\theta L^\infty \parallel  v_{-} ; L^{\infty}((0, t ], H^\lambda )\parallel \Big )  \ .
\eea

\noi The value of $\theta$ coming from Lemma 5.1 is
$$\left ( (\lambda - 1) /2 \wedge (1/8)\right ) + \mu \geq (\rho - 1)/2 \wedge (1/8)$$

\noi and we have therefore replaced it in (\ref{5.30e}) by the latter quantity, which is independent of $\lambda$. From (\ref{5.30e}) with $T_\lambda$ satisfying (\ref{5.28e}), we obtain (\ref{5.27e}) which proves the stated continuity in the interval $(0, T_\lambda )$. The extension of the continuity to the whole interval follows by an iteration argument using a simplified version of Lemma 5.1.\\

\noi \underline{Part (2)}. The proof is similar to that of Proposition 3.3. Let $A(v)$ be the linear map $v'_0 \to v'$ from $H^{\rho '}$ to ${\cal C}([0, T], H^{\rho '})$ defined by Proposition 3.4. Let $0 < \varepsilon < 3/2 - \rho$ (for instance $\varepsilon = 3/4 - \rho /2$). It follows from (\ref{3.1e}) that
\bea
\label{5.31e}
&&\parallel  A(v) v'_0 ; L^{\infty}((0, T], H^\rho  )\parallel \  \leq \overline{E} \parallel  v'_{0} ; H^\rho  \parallel  \\
&&\parallel  A(v) v'_0 ; L^{\infty}((0, T], H^{\rho + \varepsilon} )\parallel \  \leq \overline{E} \parallel  v'_{0} ; H^{\rho  + \varepsilon}\parallel 
\label{5.32e}
\eea

\noi for some constant $\overline{E} =$ Max $E(T,R)$ where the maximum is taken over $\rho ' = \rho , \rho + \varepsilon$. Let now $v_i$, $i =1,2$ be two solutions of the equation (1.21) with $v_i(0) = v_{0i} \in B(R/2)$ as obtained in Proposition 4.1. Then $v_i = A_i v_{0i}$ where $A_i = A(v_i)$. We want to show that $v_2$ tends to $v_1$ when $v_{02}$ tends to $v_{01}$ for fixed $v_{01}$. We separate again high and low frequency according to (\ref{3.25e}). We estimate
\bea
\label{5.33e}
&&\parallel  v_2 - v_1 ; L^{\infty}((0, T ], H^\rho  )\parallel \  \leq \ \parallel \left ( A_2 v_{02>} - A_1 v_{01>}\right )_> ; L^{\infty}(H^\rho  )\parallel \nn \\
&&+ \  \parallel \left ( A_2 v_{02<} - A_1 v_{01<}\right )_> ; L^{\infty} (H^\rho  )\parallel\ + \ \parallel  (v_2 - v_1)_< ; L^{\infty}(H^\rho  )\parallel  \ .
\eea

\noi Now
\bea
\label{5.34e}
&&\parallel \left ( A_2 v_{02>} - A_1 v_{01>}\right )_> ; L^{\infty} (H^\rho  )\parallel\ \leq  \  \overline{E} \sum_i \parallel  v_{0i>} ; H^\rho  \parallel \nn \\
&&\leq \  \overline{E} \left ( 2 \parallel v_{01>}  ; H^\rho  \parallel \ + \ \parallel v_{02}  - v_{01} ; H^\rho  \parallel \right )  \  .
\eea

\noi On the other hand
\bea
\label{5.35e}
&&\parallel \left ( A_2 v_{02<} - A_1 v_{01<} \right )_> ; L^{\infty} (H^\rho  )\parallel\ \leq  \  N^{- \varepsilon}  \parallel  \left ( A_2 v_{02<} - A_1 v_{01<}\right )_> ; L^{\infty} (H^{\rho  +  \varepsilon}) \parallel \nn \\
&&\leq \ N^{- \varepsilon} \sum_i \parallel A_i v_{0i<} ; L^\infty (H^{\rho  + \varepsilon})\parallel \ \leq \ \overline{E}\ N^{- \varepsilon} \sum_i \parallel  v_{0i<} ; H^{\rho  + \varepsilon}\parallel \nn \\
&&\leq \  \overline{E} \left ( 2 N^{- \varepsilon} \parallel v_{01<}  ; H^{\rho + \varepsilon} \parallel \ + \ \parallel v_{02}  - v_{01} ; H^\rho  \parallel \right )  \  .
\eea

\noi Let now $1 < \lambda < \rho$ (for instance $\lambda = (\rho + 1)/2)$. It follows from Part (1) that
\bea
\label{5.36e}
&&\parallel  (v_2 - v_1)_< ; L^{\infty}((0, T_\lambda  ], H^\rho  )\parallel \  \leq \  (N+ 1)^{\rho - \lambda }\ \parallel v_2 - v_1 ; L^\infty ((0, T_\lambda ], H^\lambda ) \parallel\nn \\
&&\leq \  (N+1)^{\rho - \lambda} C_-(\lambda ) \parallel v_{02} - v_{01} ; H^{\lambda } \parallel 
\eea

\noi for some $T_\lambda$ possibly smaller than $T$ and for some $C_- (\lambda )$ which can be read from (\ref{5.27e}). Substituting (\ref{5.34e})-(\ref{5.36e}) into (\ref{5.33e}) yields 
\bea
\label{5.37e}
&&\parallel  v_2 - v_1 ; L^{\infty}((0, T_\lambda  ], H^\rho  )\parallel \  \leq \ 2 \overline{E} \left \{  \parallel v_{01>}; H^{\rho} \parallel \ + \ N^{- \varepsilon} \parallel v_{01<}  ; H^{\rho + \varepsilon} \parallel \right \} \nn \\
&&+\  \left ( 2\overline{E} +(N+1)^{\rho - \lambda} C_-(\lambda ) \right ) \parallel v_{02} - v_{01} ; H^{\rho } \parallel \ .
\eea

\noi The two terms in the bracket tend to zero when $N$ tends to infinity by definition and by the Lebesgue dominated convergence theorem respectively, while the last term in (\ref{5.37e}) tends to zero for fixed $N$ when $v_{02}$ tends to $v_{01}$ in $H^\rho$. This completes the proof of continuity from $H^\rho$ to $L^{\infty}((0, T_\lambda  ], H^\rho  )$. The extension to the original $T$ proceeds by simpler arguments. \par\nobreak \hfill $\sq$

\noi {\bf Remark 5.1.} Using a minor generalization of Lemma 5.1, one could solve the Cauchy problem for the equation (1.21) down to $t= 0$, starting from a (sufficiently small) positive initial time.

\mysection{The Cauchy problem at time zero for u$_{\bf c}$}
\hspace*{\parindent}
In this section, we first translate the results of Sections 4 and 5 on the Cauchy problem for the equation (1.21) for $v$ into results on the Cauchy problem for the equation (1.13) for $u_c$. We then show that conversely, any sufficiently regular solution $u_c$ of (1.13) can be recovered from a suitable $v$ solution of (1.21). As an application of the latter result, we derive a uniqueness result for the equation (1.13) not making any reference to $v$. \\

\noi {\bf Proposition 6.1.} {\it Let $1 < \rho < 3/2$, let $v_0 \in H^\rho$, define $\varphi$ by (1.20) (1.14). Then there exists $T > 0$ and a unique solution $u_c \in {\cal C}((0, T], H^\rho )$ of (1.13) such that $v(t)$ defined by (1.15) satisfies the equation (1.21) and that $v \in {\cal C} ([0, T], H^\rho )$ with $v(0) = v_0$. The map $v_0 \to v$ is continuous from $H^\rho$ to $L^\infty ((0, T], H^\rho )$ and $u_c$ satisfies the estimate 
\beq
\label{6.1e}
\parallel  u_c(t) ; H^\rho \parallel\ \leq \ a_1 (1 + |\ell n\ t|)^2 
\eeq

\noi for some $a_1 \geq 0$ and for all $t \in (0, T]$.}\\

\noi {\bf Proof.} The results follow from Proposition 4.1 and 5.1 through the change of variables (1.15). The estimate (\ref{6.1e}) follows from (\ref{2.25e}).\par\nobreak \hfill $\sq$

We now derive a converse of Proposition 6.1.\\

\noi {\bf Proposition 6.2.} {\it Let $1 < \rho < 3/2$, let $I = (0, T]$ and let $u_c \in {\cal C}(I, H^\rho )$ be a solution of (1.13) in $I$ satisfying
\beq
\label{6.2e}
\parallel  u_c(t) ; H^\rho \parallel\ \leq \ a_1 (1 + |\ell n\ t|)^\alpha 
\eeq

\noi for some $a_1, \alpha \geq 0$ and for all $t \in I$. Then there exists $v \in  {\cal C} ([0, T], H^\rho )$) satisfying the equation (1.21) with $v_0 = v(0)$ such that $u_c$ is recovered from $v$ through (1.15) (1.20). Furthermore
\beq
\label{6.3e}
B(v_0) = \lim_{t\to 0} B(u_c, t) 
\eeq

\noi in $\dot{H}^\beta$ for $1/2 < \beta < 2\rho - 1/2$.}\\

\noi {\bf Proof.} We first prove the existence of the limit of $B(u_c, t)$ when $t$ tends to zero. From (\ref{2.8e}) (\ref{6.2e}), it follows that 
\beq
\label{6.4e}
\parallel  \omega^{2\sigma - 1/2} B(u_c, t) \parallel_2 \ \leq \ C\ a_1^2 (1 + |\ell n\ t|)^{2 \alpha} 
\eeq

\noi for $1/2 < \sigma \leq \rho$. From Lemma 2.7, especially (\ref{2.16e}), with $V = t^{-1} B(u_c)$, it follows that
\beq
\label{6.5e}
\parallel  \omega^{2\sigma - 5/2} \left ( B(u_c, t) - B(u_c, t_1)\right ) \parallel_2 \ \leq \ C(3 - 2 \sigma )^{-1} \ a_1^2 \ t(1 + |\ell n\ t|)^{2 \alpha} 
\eeq

\noi for $0 < t_1 \leq t$ and $1 < \sigma \leq \rho$. Interpolating between (\ref{6.4e}) and (\ref{6.5e}) yields
\beq
\label{6.6e}
\parallel  \omega^{2\sigma - 1/2- 2 \theta} \left ( B(u_c, t) - B(u_c, t_1)\right ) \parallel_2 \ \leq \ C(3 - 2 \sigma )^{-\theta} \ a_1^2 \ t^\theta (1 + |\ell n\ t| + |\ell n \ t_1|)^{2 \alpha} 
\eeq

\noi for $1 < \sigma \leq \rho$ and $0 \leq \theta \leq 1$. By the same dyadic argument as in the proof of Lemma 2.10, this implies
\beq
\label{6.7e}
\parallel  \omega^{2\sigma - 1/2- 2 \theta} \left ( B(u_c, t) - B(u_c, t_1)\right ) \parallel_2 \ \leq \ C(\sigma , \theta  ) \ a_1^2 \ t^\theta (1 + |\ell n\ t| )^{2 \alpha} 
\eeq

\noi for $0 < t_1 \leq t$, $1 < \sigma \leq \rho$ and $0 < \theta \leq 1$. From (\ref{6.4e}) (\ref{6.7e}) it follows that $B(u_c, t)$ has a limit $B_0$ in $\dot{H}^\beta$ for $1/2 < \beta < 2\rho - 1/2$ (see however Remark 2.2) and that
\beq
\label{6.8e}
\parallel  \omega^{2\sigma - 1/2- 2 \theta} \left ( B(u_c, t) - B_0)\parallel_2 \ \leq \ C(\sigma , \theta  ) \ a_1^2 \ t^\theta (1 + |\ell n\ t| )^{2 \alpha}\right .  \  . 
\eeq

\noi From (\ref{6.4e}) and (\ref{6.8e}) taken for some fixed $t$, it follows that
\beq
\label{6.9e}
\parallel  B_0 ; \dot{H}^\beta \parallel \ \leq \ C\ a_1^2 
\eeq

\noi for $1/2 < \beta < 2 \rho - 1/2$.\par

We now define $v$ by (1.15) with $\varphi = - (\ell n\ t)B_0$ so that $v$ satisfies the equation
\beq
\label{6.10e}
i \partial_t v = - t^{-1} \left ( U_\varphi^* \ B(u_c) U_\varphi - B_0 \right ) v\  .
\eeq

\noi We next prove by a variant of Propositions 3.2 and 3.3 that $v \in  {\cal C} ([0, T], H^\rho )$ and that $B(v(0)) = B_0$. We estimate $v$ by (\ref{3.5e})-(\ref{3.8e}) with $(v', \rho ')$ replaced by $(v, \rho )$. In the same way as in the proof of Proposition 3.2, using now (\ref{6.8e}) (\ref{6.9e})we obtain 
\beq
\label{6.11e}
|J_0| + |J_1| + |J_2| \leq C\ a_1^2 (1 + a_1^2 |\ell n\ t|)^4 t^\theta \left (  (1 + |\ell n \ t|)^{2 \alpha} + 1 + a_1^2 |\ell n\ t| \right ) \parallel  \omega^{\rho} v\parallel_2^2 
\eeq

\noi for $0 < \theta < \rho -1$, from which boundedness of $v$ in $L^\infty ((0, T], H^\rho )$ follows. Similarly we estimate $\parallel  v(t) - v(t_1)\parallel_2$ by (\ref{3.17e})-(\ref{3.19e}) where now 
\beq
\label{6.12e}
K_0 + K_1  \leq C\ t^{\rho /2} a_1^2  (1 + a_1^2 |\ell n \ t|)^2 (1 + |\ell n\ t|  )^{2\alpha}  \parallel  \omega^{\rho} v(t_1) \parallel_2 
\eeq

\noi which implies the continuity of $v$ in $L^2$ at $t = 0$. The regularity of $v$ then follows by the same arguments as in Propositions 3.2 and 3.3. The fact that $B_0 = B(v_0)$ is proved by the same argument as in the end of the proof of Lemma 2.10. \par \nobreak \hfill $\sq$

We finally state the uniqueness result. \\

\noi {\bf Proposition 6.3.} {\it Let $1 < \rho < 3/2$, let $I = (0, T]$ and let $u_{ci} \in {\cal C}(I, H^\rho )$, $i = 1,2$, be two solutions of the equation (1.13) in $I$ satisfying (\ref{6.2e}) for some $a_1, \alpha \geq 0$ and such that $u_{c1} - u_{c2}(t)$ tends to zero in $L^2$ when $t\to 0$. Then $u_{c1} = u_{c2}$.} \\

\noi {\bf Proof.} The main step consists in proving that $u_{c1}$ and $u_{c2}$ yield the same $B_0$. For that purpose we estimate
\begin{eqnarray*}
\parallel  \omega^{2\sigma - 1/2} \left ( B(u_{c2}, t) - B(u_{c1}, t)\right ) \parallel_2 &\leq& C \int_1^\infty d\nu\ \nu^{-2\sigma} \parallel  \omega^{2\sigma - 3/2} \Big ( |u_{c2}(t/\nu)|^2\\
&& - \  |u_{c1}(t/\nu)|^2 \Big ) \parallel_2 \end{eqnarray*}
\beq
\label{6.13e}
\leq \int d\nu\ \nu^{-2\sigma} \parallel  |u_{c2}(t/\nu)| -  |u_{c1}(t/\nu)|\parallel_2\left ( \parallel  u_{c1}(t/\nu)\parallel_r\ + \ \parallel  u_{c2}(t/\nu)\parallel_r \right ) \nn \\
\eeq

\noi for $1/2 < \sigma \leq 3/4$, with $3/2 - 3/r = 2\sigma$. We next estimate 
$$|u_c| \leq |(U-1)\ e^{-i\varphi} v| + |v|$$

\noi for $u_c = u_{ci}$, $i = 1,2$, so that 
\bea
\label{6.14e}
\parallel u_c \parallel_r &\leq &C\left ( t^\mu \parallel  \omega^{2(\sigma + \mu )} e^{-i\varphi} v  \parallel_2 \ + \ \parallel \omega^{2\sigma} v\parallel_2 \right )\nn \\
&\leq &C\left ( t^\mu \parallel  e^{-i\varphi} ; M^{2(\sigma + \mu )} \parallel \ \parallel \omega^{2(\sigma + \mu )} v  \parallel_2 \ + \ \parallel \omega^{2\sigma} v\parallel_2 \right )\nn \\
&\leq &C\left ( 1 + t^\mu (1 + a_1^2|\ell n\ t|)^2 \right ) \parallel v ; L^\infty (I, H^\rho )  \parallel                                                                                                                                                                                                                                                                                                                               
\eea

\noi for $\sigma + \mu \leq \rho /2$, which can be achieved with $\mu > 0$ for $\sigma < \rho /2$. It follows from  (\ref{6.13e})  (\ref{6.14e}) through the Lebesgue dominated convergence theorem that $B(u_{c2}, t) - B(u_{c1}, t)$ tends to zero in $\dot{H}^\beta$ for $1/2 < \beta < \rho - 1/2$ since $|u_{c2}| - |u_{c1}|$ tends to zero in $L^2$ when $t$ tends to zero. This implies that $B_{02} = B_{01}$ so that $\varphi_2 = \varphi_1$ and $v_i$, $i=1,2$, satisfy the same equation (1.21). \par

On the other hand $\parallel v_2(t) - v_1(t)\parallel_2 \ =\  \parallel u_{c2}(t) -  u_{c1}(t)\parallel_2$ tends to zero by assumption, so that $v_{02} = v_{01}$. The uniqueness result of Proposition 4.1 then implies that $v_2 = v_1$ and therefore $u_{c2} = u_{c1}$.\par\nobreak \hfill $\sq$


\begin{thebibliography}{99}
\bibitem{1r} J. Bergh, J. L\"ofstr\"om, Interpolation Spaces, Springer, Berlin, 1976.

\bibitem{2r} J. Ginibre, G. Velo, Long range scattering and modified wave operators for some Hartree-type equations I,  Rev. Math. Phys. {\bf 12} (2000) 361-429.

\bibitem{3r} J. Ginibre, G. Velo, Long range scattering and modified wave operators for the Wave Schr\"odinger system, Ann. H. P. {\bf 3} (2002) 537-612. Id II, Ann. H. P. {\bf 4} (2003) 973-999. Id III, Dynamics of PDE, {\bf 2} (2005) 101-125. 

\bibitem{4r} K. Nakanishi, Modified wave operators for the Hartree equation with data, image and convergence in the same space, Commun. Pure Appl. Anal. {\bf 1} (2002) 237-252.

\bibitem{5r} K. Nakanishi, Modified wave operators for the Hartree equation with data, image and convergence in the same space II, Ann. H.P. {\bf 3} (2002) 503-535.

\bibitem{6r} W. Strauss, Nonlinear Wave Equations, CMBS Lecture Notes 73, Am. Math. Soc., Providence, 1989.
\end{thebibliography}
\end{document}